\documentclass[12pt]{article}

\usepackage{amssymb,amsmath,amsthm,amstext} 

\newcommand{\nek}{\newcommand}
\nek{\renek}{\renewcommand}

\renek{\thesubsection}{\arabic{subsection}}

\nek{\parf}[1]{\subsection{\boldmath#1}}
\nek{\punk}[1]{\subsubsection{\boldmath#1}}

\DeclareMathAlphabet{\cur}{U}{eur}{m}{n}
\DeclareMathAlphabet{\skr}{U}{eus}{m}{n}
\SetMathAlphabet{\cur}{bold}{U}{eur}{b}{n}
\SetMathAlphabet{\skr}{bold}{U}{eus}{b}{n}
\DeclareMathAlphabet{\bma}{OML}{cmm}{b}{it}

\nek{\bfit}{\bfseries\itshape}
\nek{\bfsl}{\bfseries\slshape}

\nek{\vyk} [1] {}

\nek{\itsep}{\itemsep=0.4ex plus 0.15ex minus 0.15ex}
\nek{\tenu}[1]{
\def\theenumi{#1}
\itsep
}
\newtheorem{theore}             {Theorem} 
\newtheorem{corollar}  [theore]{Corollary}
\newtheorem{propo}     [theore]{Proposition}
\newtheorem{lemm}      [theore]{Lemma}
\newtheorem{cla}       [theore]     {Claim} 

\theoremstyle{definition}
\newtheorem{defn}      [theore]{Definition}
\newtheorem{rem}       [theore]{Remark}

\newtheorem*{prF}{Proof}               

\newtheorem{que}  {Question}
\nek{\bpro}{\begin{propo}}
\nek{\epro}{\end{propo}}
\nek{\bcl} {\begin{cla}}
\nek{\ecl} {\end{cla}}
\nek{\bcor}{\begin{corollar}}
\nek{\ecor}{\end{corollar}}
\nek{\bdf} {\begin{defn}}
\nek{\eDf} {\end{defn}}
\nek{\edf} {\qed\end{defn}}
\nek{\edF} {\end{defn}}
\nek{\ble} {\begin{lemm}}
\nek{\ele} {\end{lemm}}
\nek{\bte} {\begin{theore}}
\nek{\ete} {\end{theore}}
\nek{\bre} {\begin{rem}} 
\nek{\ere} {\qeD\end{rem}} 
\nek{\bqe} {\begin{que}} 
\nek{\eqe} {\qed\end{que}} 
\nek{\bpf} {\begin{prF}} 
\nek{\epf} {\qed\end{prF}} 
\nek{\ePf} {\end{prF}} 
\nek{\qeDD} [1] 
{\hfill\text{$\mtho\qed$~({\sl #1\/}\hspace{0.1ex})}}
\nek{\epF} [1] {\qeDD{#1}\end{prF}} 

\nek{\bde}{\begin{description}}
\nek{\ede}{\end{description}}
\nek{\ben}{\begin{enumerate}}
\nek{\een}{\end{enumerate}}
\nek{\bit}{\begin{itemize}}
\nek{\eit}{\end{itemize}}
\nek{\bay}{\begin{array}}
\nek{\eay}{\end{array}}
\nek{\fF} {{\bf F}}
\nek{\fG} {{\bf G}}
\nek{\Gd} {\fG_\fda}
\nek{\Fs} {\fF_\fsg}
\nek{\Gds} {\fG_{\fda\fsg}}
\nek{\Fsd} {\fF_{\fsg\fda}}

\nek{\fa} {{\bf a}}
\nek{\fb} {{\bf b}}

\nek{\ZFC}{{\bf ZFC}}
\nek{\zhc}{{\bf ZFHC}}
\nek{\iesp}{\hspace{0.3ex}}

\nek{\ie} {\hbox{\sl i.\iesp e.}}
\nek{\eg} {\hbox{\sl e.\iesp g.}}
\nek{\lsc} {\hbox{l.\iesp s.\iesp c.}}
\nek{\er} {{ER}}
\nek{\vrt} {\hbox{w.\iesp r.\iesp t.}}
\nek{\bm} {{BM}}
\nek{\etc} {{\sl etc}}
\nek{\dd}[1]{$\mtho\hspace*{0.2ex}{#1}$-}

\nek{\ran}  {\mathop{\tt ran}}
\nek{\dom}  {\mathop{\tt dom}}
\nek{\tsup} {\mathop{\tt sup}}
\nek{\tinf} {\mathop{\tt inf}}
\nek{\Ord}  {{\tt{Ord}}}
\nek{\Exh}  {{\tt{Exh}}}
\nek{\Mod}  {\mathop{\tt{Mod}}}
\nek{\tlim} {\mathop{\tt lim}}
\nek{\tlis} {\mathop{\tt lim\hspace{0.3ex}sup}}
\nek{\card} {\mathop{\tt card}}
\nek{\otp}  {\mathop{\tt otp}}
\nek{\lh}   {\mathop{\tt lh}}

\nek{\ccc}{{\sc ccc}}
\nek{\al} {\alpha}
\nek{\ba} {\beta}
\nek{\ga} {\gamma}
\nek{\Da} {\Delta}
\nek{\la} {\lambda}
\nek{\La}{\Lambda}
\nek{\ve} {\varepsilon}
\nek{\vt} {\vartheta}
\nek{\vpi}{\varphi}
\nek{\vpy}{\vpi_\iy}
\nek{\vT} {\Theta}
\nek{\om} {\omega}
\nek{\omi} {\om_1}
\nek{\omm} [1] {\om^{\om^{#1}}}
\nek{\lom}{^{<\om}}
\nek{\sd}   {\mathbin{\Da}}
\nek{\alo} {{\aleph_0}}

\nek{\fs}[2]{{\bf\iSg}^{#1}_{#2}}
\nek{\fp}[2]{{\bf\iPi}^{#1}_{#2}}
\nek{\fd}[2]{{\bf\iDa}^{#1}_{#2}}
\nek{\is}[2]{\iSg^{#1}_{#2}}
\nek{\ip}[2]{\iPi^{#1}_{#2}}
\nek{\id}[2]{\iDa^{#1}_{#2}}
\nek{\iSg}{{\mathchar"7106}}
\nek{\iPi}{{\mathchar"7105}}
\nek{\iDa}{{\mathchar"7101}}
\mathchardef\deltA ="710E
\mathchardef\sigmA ="711B
\mathchardef\xI ="7118
\mathchardef\taU ="711C

\nek{\fsg} {{\cur{\sigmA}}}
\nek{\fda} {{\cur{\deltA}}}

\nek{\dds}{\dd\fsg}

\nek{\BBB}{\hspace{0.05ex}}
\nek{\dC}{{\BBB{\mathbb C}\BBB}}
\nek{\dG}{{\BBB{\mathbb G}\BBB}}
\nek{\dN}{{\BBB{\mathbb N}\BBB}}
\nek{\dP}{{\BBB{\mathbb P}\BBB}}
\nek{\dS}{{\BBB{\mathbb S}\BBB}}
\nek{\dZ}{{\BBB{\mathbb Z}\BBB}}
\nek{\dX}{{\BBB{\mathbb X}\BBB}}
\nek{\dY}{{\BBB{\mathbb Y}\BBB}}
\nek{\dV}{{\BBB{\mathbb V}\BBB}}
\nek{\dn}{2^\dN}
\nek{\dvp} {\dV^+}

\nek{\cI} {{\cal I}} 
\nek{\cZ} {{\cal Z}}
\nek{\cJ} {{\cal J}} 
\nek{\cO} {{\cal O}} 
\nek{\cL}{{\BBB{\cal L}\BBB}}
\nek{\cP}{{\BBB{\skr P}\BBB}}
\nek{\cW}{{\BBB{\skr W}\BBB}}
\nek{\cU}{{\BBB{\skr U}\BBB}}
\nek{\cG}{{\BBB{\skr G}\BBB}}
\nek{\pn}{\cP(\dN)}
\nek{\pws}  [1] {\cP(#1)}
\nek{\mm}{{\BBB{\mathfrak M}\BBB}}
\nek{\Zo}  {{\cZ_0}}

\nek{\iz} [2] {\cI_{#1}^{#2}}
\nek{\iw} [2] {\cW_{#1}^{#2}}

\nek{\ir} [2] {{[#1,#2)}}
\nek{\iry}[1] {{[#1,\iy)}}
\nek{\opl} {\oplus}
\nek{\ap}  {\cdot}
\nek{\app} {{\hspace{0.2ex}{\cdot}\hspace{0.2ex}}}
\nek{\dm}  {$$}
\nek{\sus} {\mathopen{\exists\hspace{0.35ex}}}
\nek{\kaz} {\mathopen{\forall\hspace{0.35ex}}}
\nek{\imp} {\Longrightarrow} 
\nek{\eqv} {\Longleftrightarrow} 
\nek{\ti}  {\times} 
\nek{\mo}  {\models} 
\nek{\sq}  {\subseteq}
\nek{\su}  {\subset}
\nek{\sneq}{\subsetneqq}
\nek{\we}  {{\mathbin{\hspace*{0.2ex}^\wedge}}}
\nek{\obr} {^{-1}}
\nek{\dif} {\setminus}
\nek{\res} {\mathop{\restriction}}
\nek{\pu}  {\emptyset}
\nek{\iy}  {\infty}
\nek{\piy} {+\iy}
\nek{\nin} {\not\in}
\nek{\limp}{\,\imp\,}
\nek{\leqv}{\,\eqv\,}
\nek{\onto}{\stackrel{{\rm onto}}{\longrightarrow}}
\nek{\ang} [1] {\langle #1\rangle}
\nek{\stk} [2] {\ang{#1\hspace{0.3ex};\hspace{0.1ex}#2}}
\nek{\sis} [2] {\ans{#1}_{#2}}
\nek{\ans} [1] {\{\hspace{0.01ex}#1\hspace{0.01ex}\}} 
\nek{\itla} {\item\label}

\nek{\seq}[2] {(#1)_{#2}}

\nek{\ifi}  {{\cur{Fin}\hspace*{0.1ex}}}
\nek{\frt}  {{\cur{Fr}\hspace*{0.1ex}}}
\nek{\fio} {\ifi\ti0}
\nek{\ofi} {0\ti\ifi}

\nek{\skl} {\hbox{\mtho\large$($}}
\nek{\skp} {\hbox{\mtho\large$)$}}

\nek{\rD}  {\mathbin{\sf D}}
\nek{\rF}  {\mathbin{\sf F}}
\nek{\rE}  {\mathbin{\sf E}}
\nek{\rP}  {\mathbin{\sf P}}
\nek{\rR}  {\mathbin{\sf R}}
\nek{\rT}  {\mathbin{\sf T}}
\nek{\rtd} {\rT_2}
\nek{\rei} {\rE_{\cI}}
\nek{\rej} {\rE_{\cJ}}

\nek{\nE}  {\mathbin{{\not\hspace{-0.35ex}\sf E}}}
\nek{\nF}  {\mathbin{{\not\hspace{-0.25ex}\sf F}}}

\nek{\Eo}  {\rE_0}
\nek{\reo} {\rE_{\hspace{-1.0pt}\Zo}}

\nek{\dde} {\dd{\rE}}
\nek{\ddf} {\dd{\rF}}
\nek{\ef} {\dd{\rE,\rF}}

\nek{\reb} {\le_{\rm B}}
\nek{\eqb} {\approx_{\rm B}}
\nek{\rez} {\rE_{\hspace{-1.0pt}\cZ}}

\nek{\fps} [3] {{\prod_{#3}{#2}\,/\,{#1}}}
\nek{\fpd} [2] {{#1}\otimes{#2}}
\nek{\rkb} [1] {|#1|_{\rm CB}}

\nek{\spa} {\dS}
\nek{\spn} {\spa^\dN}

\nek{\poq}{\underline}
\nek{\nad}{\overline}

\nek{\tsc}[1]{{\hbox{\footnotesize\sc{#1}}}}

\nek{\aprb}{\approx_{\tsc b}}
\nek{\ismb}{\cong_{\tsc b}}

\nek{\rebs} {<_{\tsc b}}
\renek{\reb} {\le_{\tsc b}}
\renek{\eqb} {\sim_{\tsc b}}

\nek{\Ii}{\cI_1}
\nek{\It}{\cI_3}
\nek{\Ei}{\rE_1}
\nek{\Et}{\rE_3}

\nek{\ong} {1_\dG}
\nek{\lo} [3] {\cO(#3,#1,#2)}
\nek{\ler}[2] {\mathbin{\sim^{#1}_{#2}}}
\nek{\aer}[2] {\mathbin{\rE_{#1}^{#2}}}
\nek{\ergx}{\aer\dG\dX}
\nek{\egx} {\ergx}
\nek{\isg} {{S_\iy}}
\nek{\dnn}{\dN^\dN}
\nek{\rr} [2] {\rR^{#1}_{#2}}
\nek{\sym}{\ler} 

\nek{\rav}[1] {\rD(#1)}
\nek{\ek}[2] {[#1]_{{#2}}}
\nek{\eke}[1] {\ek{#1}{\rE}}

\nek{\ur}{_{\tt right}}
\nek{\ul}{_{\tt left}}

\nek{\cli} {CLI} 
\nek{\di} [1] {{#1}^\#}
\nek{\drE} {\mathbin{\di\rE}}
\nek{\drF} {\mathbin{\di\rF}}

\nek{\mox} {\moq} 
\nek{\loa} [1] {j_{#1}}
\nek{\ism} [1] {\cong_{#1}}
\nek{\moq} [1] {\Mod_{#1}}
\nek{\hfn} {{\rm HF}(\dN)}
\nek{\tce} [1] {{\rm TC}_\ve(#1)}
\nek{\ihf} {\simeq_{\hfn}}
\nek{\symr}{\equiv}
\nek{\rrt} [3] {\symr_{#2#3}^{#1}}
\nek{\rrq} [5] {{#4}\symr_{#2#3}^{#1}{#5}}
\nek{\rrQ} [5] {{#4}\symr_{#2\,,\,#3}^{#1}{#5}}
\nek{\nrq} [5] {{#4}\not\symr_{#2#3}^{#1}{#5}}

\nek{\rro} [5] {\ang{#2,#4}\symr^{#1}\ang{#3,#5}}
\nek{\rrO} [1] {\symr^{#1}}
\nek{\lww} {\cL_{\omi\om}}

\nek{\forc}
{\mathrel{{|\hspace{-0.2ex}|\hspace{-1.1ex}-\hspace{-1.7ex}-}}}

\newlength{\dxii}
\nek{\fC} {{\bf C}}
\nek{\pg} {\fC_\dG}
\nek{\px} {\fC_\dX}

\nek{\gen} {gen.\ }
\nek{\genp}{gen.}
\nek{\hg}  {h.\hspace{0.4ex}gen.\ }
\nek{\hgp} {h.\hspace{0.4ex}gen.}

\nek{\rfy} {\rF^\iy}

\nek{\incl} [1] {\mathop{\text{\sc Int}}\overline{#1}}  

\nek{\PP} {pinned}
\nek{\PPP}{Pinned}

\nek{\zO} {,\linebreak[0]}
\nek{\zi} {,\linebreak[0]\,}
\nek{\zd} {,\linebreak[0]\:}
\nek{\zT} {,\linebreak[0]\;}
\nek{\zq} {,\linebreak[0]\;\,}

\nek{\aci}[1]{\addtocounter{enumi}#1}

\nek{\dPP} {\dP\ti\dP}
\nek{\dSP} {\dS\ti\dP}
\nek{\xib} {{\bma\xI\hspace{0.1ex}}}
\nek{\sgb} {\bma\sigmA}
\nek{\ddp} {\dd\dP}
\nek{\coll} {\text{\sc Coll}}

\nek{\dCP} {\dC\ti\dP}
\nek{\dCC} {\dC\ti\dC}

\nek{\tal}{\xib\ul}
\nek{\tar}{\xib\ur}

\nek{\dNN}{(\dn){\vphantom{\dN}}^\dN}

\nek{\sgbl}{\sgb\ul}
\nek{\sgbr}{\sgb\ur}

\nek{\mtho}{\mathsurround=0mm}
\nek{\msur}{\hspace{-1\mathsurround}}
\nek{\psur}{\hspace{0.3\mathsurround}}
\nek{\dsur}{\hspace{-0.3\mathsurround}}
\nek{\hsur}{\hspace{-0.5\mathsurround}}
\nek{\noi}{\noindent}
\nek{\vom}{\vspace{1mm}}
\nek{\vtm}{\vspace{2mm}}

\begin{document}

\title{Some new results on Borel irreducibility
of equivalence relations~\thanks
{\ The results of this paper were partially presented at 
LC`1999 (Utrecht), LC`2000 (Paris), Lumini Meeting in 
set theory (2000), Caltech and UCLA seminars, and 
P.~S.~Novikov's Memorial Conference (Moscow, 2001).}}

\date{March 2002}

\author{Vladimir Kanovei~\thanks
{\ 
Moscow, {\tt kanovei@math.uni-wuppertal.de}, 
{\tt kanovei@mccme.ru}.
Suported by grants of  
{\rm DFG 101/10-1, NSF DMS 96-19880},   
universities of Bonn, Wuppertal, and Caltech.
}
\and
Michael Reeken~\thanks
{\ University of Wuppertal, Germany, 
{\tt reeken@math.uni-wuppertal.de}.}
}

\maketitle

\begin{abstract}
We prove that orbit equivalence relations 
(\er s, for brevity) of 
generically turbulent Polish 
actions are not Borel reducible to \er s of a family which 
includes Polish actions of  
$\isg,$ the group of all permutations of $\dN,$ 
and is closed under the Fubini product modulo the 
ideal $\ifi$ of all finite sets, and some other operations.  
Our second main result shows that $\rtd,$ an 
equivalence relation called 
``the equality of countable sets of the reals'', 
is not Borel reducible to another family of \er s which 
includes continuous actions of Polish \cli\ groups, 
Borel equivalence relations with $\Gds$ classes, some 
ideals, and is closed under the Fubini product over $\ifi.$ 
Both results and their corollaries extend some earlier 
irreducibility theorems by Hjorth and Kechris.  
\end{abstract}


\subsection*{Introduction}
\label i

Classification problems for different types of mathematical 
structures have been in the center of interests in 
descriptive set theory since the beginning of the 90s. 
Suppose that $X$ is a class of mathematical structures, 
identified modulo an equivalence relation $\rE.$ 
This can be, \eg, countable groups modulo  
the isomorphism relation, or unitary operators over a 
fixed space $\dC^n$ modulo conjugacy,   
or probability measures over a fixed Polish space modulo 
the identification of measures having the same null sets, 
or, for instance, reals modulo the Turing reducibility.  
(The examples are taken from Hjorth's book~\cite h and 
Kechris' survey paper~\cite{ndir}, where many more examples 
are given.) 
Suppose that $Y$ is another class of mathematical structures, 
identified modulo an equivalence relation $\rF.$ 
The classification problem is then to find out whether 
there exists a {\it definable\/}, or {\it effective\/} 
injection $\vT:{X/\rE}\to {Y/\rF}.$ 
Such a map $\vT$ can be seen as a classification of 
objects in $X$ in terms of objects in $Y,$ in a way which 
respects quotients over $\rE$ and $\rF.$  
Its existence can be a result of high importance, for instance 
when objects in $Y$ are of mich simpler mathematical nature 
than those in $X$.

In many cases, it turns out that the classes of structures 
$X$ and $Y$ can be considered as Polish 
(\ie, separable complete metric) 
spaces, so that $\rE,\,\rF$ become Borel 
or, more generally, analytic (as sets of pairs) relations,  
while reduction maps are usually required 
to be Borel~\footnote 
{\ 
That is, with Borel graphs. 
Baire measurable maps and  
reductions satisfying certain algebraic requirements are 
also considered \cite{aq}, as well as $\fd12$ and more 
complicated reductions \cite{k:sol,k:ulm}, 
however they are not in the scope of this paper.}. 
In this case the problem can be studied by methods of  
descriptive set theory, where it takes the form: 
if $\rE,\,\rF$ are Borel or analytic equivalence relations 
on Polish spaces, resp., $X,\,Y,$ does there exist a Borel 
{\it reduction\/} 
$\rE$ to $\rF,$ \ie, a Borel map $\vt:X\to Y$ satisfying 
${x\rE x'}\eqv{\vt(x)\rF\vt(x')}$ 
for all $x,\,x'\in X.$ 
If such a map $\vt$ exists then $\rE$ is said to be 
{\it Borel reducible\/} to $\rF.$
Studies on Borel and analytic equivalence relations 
(\er s, for brevity) under Borel reducibility  
by methods of descriptive set theory, revealed 
a remarkable structure 
of reducibility and irreducibility theorems between 
\er s of different types 
(we can cite \cite{aq,h:orb,h,hk:nd,ndir} as a partial 
account of the results obtained). 
Our paper belongs to this research direction. 

Our first main theorem establishes Borel irreducibility 
between two large classes of \er s. 
{\sl Class 1\/} consists of \er s induced by  
{\it generically turbulent\/} Polish actions~\footnote
{\label{pol}\ 
That is, comeager many orbits, 
and even local orbits, are somewhere dense, see  
Definition~\ref{df:durb}.}.  
Hjorth~\cite h proved that no \er\ of this class 
is Borel reducible to an \er\ of {\sl Class~2\/} 
which consists of orbit \er s of Polish actions of 
$\isg,$ the group of all permutations of $\dN.$ 
(This result is also known in the form: generically turbulent 
\er s are not classifiable by countable structures, see 
comments in \ref{loac}.)

One of possible proofs of Hjorth's theorem is as follows. 
{\it First\/}, any \er, Borel reducible to an \er\ in Class~2, 
is then Borel reducible, at least on a comeager set, to an   
\er\ in {\sl Class~3\/}~\footnote
{\ Introduced essentially by H.~Friedman~\cite{fs,frid}.}, 
which consists of those \er s that can be obtained from 
equalities on Polish spaces using the operation of 
countable power $\rE^\iy.$ 
{\it Second\/} (this involves Hjorth's turbulence theory), 
no \er\ in Class~1 is Borel reducible to a \er\ in Class~3, 
even on a comeager set. 
Our Theorem~\ref{m1} generalizes the second part. 
We consider {\sl Class~4\/}, containing all \er s which 
can be obtained from equalities on Polish spaces with the 
operations of 
$1)\msur$ countable union (if it results in a \er) 
of \er s in the same space, 
$2)\msur$ Fubini product 
$\fps{\ifi}{\rE_{k}}{k\in\dN}$ modulo the ideal 
$\ifi$ of all finite subsets of $\dN,$ and 
$3)\msur$ the countable power $\rE^\iy$ 
(see \ref{new:er} for exact definitions). 
Class~4 includes Class~3, of course, but contains many more 
various \er s, especially those defined using Fubini 
products, for instance, 
all \er s induced by generalized Fr\'echet, 
indecomposable, and Weiss ideals (see \ref{new:er}).

\bte
\label{m1}
\er s in Class 1 are not Borel reducible 
(even not reducible by Baire measurable functions) 
to \er s in Class 4. 
\ete

The proof (Section~\ref{thm1}) involves the induction on 
the construction of \er s in Class~4 with the help of 
the operations indicated. 
The technique of the turbulence theory will play the key 
role in the proof, in particular, the key step will be to 
prove that any \er\ in Class~1 is generically ergodic 
\vrt\ any \er\ in Class~4 (Theorem~\ref{Turb}). 
As an application of this result, we derive the 
abovementioned Hjorth's theorem in a few rather simple 
steps in Section~\ref{appl}.

Amongst the habitants of Class 1 we have \er s of the 
form $x\rei y$ iff 
${x\sd y}\in\cI,$ for all $x,\,y\in\pn,$ 
where $\cI$ is an ideal on $\dN.$  
Any ideal $\cI\sq\pn$ is obviously an Abelian group 
with the symmetric difference $\sd$ as the group 
operation, and $\rei$ is induced by the 
shift action of $\cI$ on $\pn$ by $\sd.$ 
Kechris~\cite{rig} proved that this action is turbulent 
provided $\cI$ is a Borel P-ideal~\footnote
{\label{pideal}\ 
An ideal $\cI$ on $\dN$ is a {\it P-ideal\/} if for any 
sequence of sets $x_n\in\cI$ there is a set $x\in\cI$ such 
that $x_n\dif x$ is finite for any $n.$ 
For Borel ideals, this is equivalent to {\it polishability\/}, 
\ie, the existence of a Polish topology on $\cI$ which 
converts $\stk\cI\sd$ in a Polish group.}~, 
with few exceptions mentioned below. 
This allows us to prove the following theorem in \ref{thm2} 
as a corollary of Theorem~\ref{m1}. 

\bte
\label{m2}
If\/ $\cZ$ is a non-trivial~\footnote
{\label{nontr}\ 
Here, containing all singletons $\ans n\zT n\in\dN,$ 
and different from $\pn$.} 
Borel P-ideal on $\dN$ then\/ $\rez$ is not Borel 
reducible to a \er\ in Class~4 unless\/ 
$\cZ$ is $\ifi$ or a trivial variation of\/ $\ifi,$ 
or\/ $\cZ$ is isomorphic to\/ $\cI_3=\ofi$ via a bijection 
between the underlying sets. 
\ete

Borel P-ideals form a widely studied class,  
which includes, for instance, $\ifi,$ the ideal
$\cI_3=\ofi$ of all sets $x\sq\dN^2$ such that every 
cross-section $(x)_n=\ans{k:\ang{n,k}\in x}$ is finite, 
and {\it trivial variations of\/} $\ifi,$ \ie, ideals 
of the form $\cI_W=\ans{x\sq\dN:x\cap W\in\ifi},$ 
where $W\sq\dN$ is infinite (see \ref{new:i}), 
as well as summable ideals, density ideals, and many more 
(see \cite{aq,sol'}).
Easily the \er s $\Eo=\rE_\ifi$ and $\Et=\rE_{\It}$  
induced by the ideals $\ifi$ and $\cI_3,$ 
and those induced by trivial variations 
of $\ifi,$ belong to Class~4, 
thus, the exclusion of $\ifi\zT\It,$ and 
trivial variations of $\ifi$ in Theorem~\ref{m2} is 
necessary and fully motivated. 

Note that a weaker form of Theorem~\ref{m2}, 
with Class~3 instead of Class~4, is implicitly contained 
in Kechris~\cite{rig}. 
A very partial result, that $\reo,$ the \er\ associated 
with the density-0 ideal $\Zo,$ is not Baire reducible to 
any \er\ of Class~3, was announced in \cite{fs}. 
(Friedman gives a proof in \cite{frid}.)  

The final Section~\ref{rtd} is written in attempt to 
obtain results of the opposite character, \ie, that \er s 
in Class~4 are not Borel reducible to, say, turbulent or 
some other \er s of different nature. 
This is a comparably less developed area, and perhaps only 
one special theorem of this sort is known: 
Hjorth~\cite{h:orb} proved that $\rtd,$ the \er\ defined on 
countable sequences of the reals so that 
${\sis{x_n}{}}\rtd{\sis{y_n}{}}$ iff 
$\ans{x_n:n\in\dN}=\ans{y_n:n\in\dN}$~\footnote 
{\ $\rtd,$ sometimes denoted $\rF_2,$ as in \cite h, 
is often called ``the equality of countable sets of reals''; 
it belongs to Class~3 and is one of the most important 
Borel \er s.}
, 
is not Borel reducible to any \er\ induced by a continuous 
action of a Polish group which admits a compatible complete 
left-invariant metric (a {\sl\cli\/} group; 
this includes, for instance, Polish {\it Abelian\/} groups). 
It can be expected that $\rtd$ is not Borel reducible to 
any Borel action of a Borel Abelian group, but this is 
still open, even \vrt\ shift \dd\sd actions of 
Borel ideals. 

A possible way to solve the problem is connected with the 
following condition of a \er\ $\rE$ 
(implicitly in \cite{h:orb}): 
{\it for any forcing notion\/ $\dP$ and any \ddp term $\xib,$ 
if\/ $\dPP$ forces\/ 
$\tal\rE\tar$ then there is a real\/ $x$ in the ground 
universe such that\/ $\dP$ forces\/ $x\rE\xib.$\/} 
We call {\it\PP\/} all \er s satisfying this condition.  
Note that $\rtd$ is not \PP and not Borel reducible to 
any analytic \PP\ \er.  
We prove in Section~\ref{rtd} that 
$1)\msur$ \er s induced by Polish actions of \cli\ 
groups are \PP\  
(our proof is a simplification of Hjorth's proof in~\cite{h:orb}), 
$2)\msur$ 
Borel \er s whose equivalence classes are $\Gds$ sets are 
\PP\ (based on an idea extended to us by Hjorth), 
$3)\msur$ 
\er s associated with exhaustive ideals of sequences 
of submeasures on $\dN$ (not all of them are Polishable) 
are \PP, 
$4)\msur$ 
Fubini products of analytic \PP\ \er s modulo $\ifi$ are 
\PP. 
All of those \er s do not Borel reduce $\rtd.$ 
For instance all \er s induced by Fr\'echet ideals are 
\PP\ and do not Borel reduce $\rtd.$

\parf{Preliminaries}
\label{not}

This Section contains a review of basic notation 
involved in the formulations and proofs of our main 
results.

\punk{Descriptive set theory}
\label{DST}

Some degree of knowledge of the theory of Borel and analytic 
sets 
in Polish (\ie, complete separable metric) spaces is assumed. 
Recall that {\it analytic\/} sets  
(also known as Suslin, A-sets, or $\fs11$)  
are continuous images of Borel sets. 
Any Borel set is analytic, but 
(in uncountable Polish spaces) not conversely. 

A map $f$ (between Borel sets in Polish spaces) is 
called {\it Borel\/} iff its graph is a Borel set, or, 
which is the same, if all 
\dd fpreimages of open sets are Borel. 
{\it BM\/} always means {\it Baire measurable\/}. 
A map $f:\dX\to \dY$ is BM iff all \dd fpreimages of 
open sets in $\dY$ have the Baire property in $\dX$ 
(\ie, are equal to open sets modulo meager sets, that 
is, sets of the 1st category). 
Any such a map is continuous on a dense $\Gd$ set 
$D\sq\dX$ ($\dX,\,\dY$ are supposed to be Polish). 

Superpositions of Borel maps are easily Borel. 
Generally speaking, this is not true for BM 
maps, however, we have a useful partial result.

\ble
\label{BBM}
If\/ $\dX,\,\dY,\,\dZ$ are Polish spaces, 
$f:\dX\to \dY$ is BM, 
and\/ $g:\dY\to \dZ$ is Borel, then the superposition\/ 
$f\circ g:\dX\to\dZ$ is BM.
\ele
\bpf
By definition, \dd gpreimages of open subsets of $\dZ$ 
are Borel in $\dY,$ whose \dd fpreimages are  
Borel combinations of sets having the Baire property.
\epf

\punk{Equivalence relations}
\label{new:er}

\er\ means: equivalence relation. 
By $\rav X$ we denote the equality on $X$ considered as a \er.

Let $\rE$ be a \er\ on a set $X.$ 
Then $\eke y=\ans{x\in X:y\rE x}$ is 
{\it the\/ \dde class\/} of any element $y\in X.$   
A set $Y\sq X$ is 
{\it pairwise \dde equivalent\/}, 
if $x\rE y$ holds for all $x,\,y\in Y$.   

Suppose that $\rE$ and $\rF$ are \er s on Polish spaces 
resp.\ $\dX,\,\dY.$  
\bit
\item[$\ast$]\msur
$\rE\reb\rF$ ({\it Borel reducibility\/}; sometimes 
they write ${\dX/\rE}\reb{\dY/\rF}$)
means that there is a Borel map 
$\vt:\dX\to\dY$ (called {\it reduction\/}) 
such that ${x\rE y}\eqv {\vt(x)\rF\vt(y)}$ for all 
$x,\,y\in\dX$;

\item[$\ast$]\msur
$\rE\eqb\rF$ means that $\rE\reb\rF$ and $\rF\reb\rE$ 
({\it Borel bi-reducibility\/});  

\item[$\ast$]\msur
$\rE\rebs\rF$ means that $\rE\reb\rF$ but $\rF\not\reb\rE$ 
({\it strict Borel reducibility\/}).  
\eit

The following operations over \er s on Polish 
spaces are considered.

\ben
\tenu{(e\arabic{enumi})} 
\itla{oe:beg}
\label{cun}
{\it the countable union\/} (if it results in a \er) and 
{\it the countable intersection\/} of \er s on one and the 
same space; 

\itla{cdun}
{\it the countable disjoint union\/} 
$\rF=\bigvee_{k\in\dN}\rF_k$ 
of \er s $\rF_k$ on Polish spaces $\spa_k$ is a \er\ 
on $\spa=\bigcup_k(\ans k\ti\spa_k)$ 
(with the Polish topology generated by sets of the 
form $\ans k\ti U,$ where $U\sq \spa_k$ is open) 
defined as follows: $\ang{k,x}\rF \ang{l,y}$ iff 
$k=l$ and $x\rF_k y$~\footnote 
{\ If $\spa_k$ are pairwise disjoint and open in 
$\spa'=\bigcup_k\spa_k$ then we can equivalently 
define $\rF=\bigvee_k\rF_k$ on $\spa'$ so that 
$x\rF y$ iff $x,\,y$ belong to the same $\spa_k$ 
and $x\rF_k y$.} ; 

\itla{prod} 
{\it the product\/} $\rP=\prod_k\rF_k$ 
of \er s $\rF_k$ on spaces $\spa_k$ is a \er\ on 
$\prod_k\spa_k$ defined so that  
$x\rP y$ iff $x_k\rF_k y_k$ for all $k,$ 
in particular, if $\rE,\,\rF$ live in, resp., $\dX,\,\dY$ 
then $\rP=\rE\ti\rF$ is defined on $\dX\ti\dY$ so that  
$\ang{x,y}\rP\ang{x',y'}$ iff $x\rE x'$ and $y\rF y'$; 

\itla{fp}
{\it the Fubini product\/} 
$\fps{\cI}{\rF_{k}}{k\in\dN}$ 
of \er s $\rF_k$ on spaces $\spa_k,$ 
modulo an ideal $\cI$ on $\dN,$ 
is a \er\ on the product space 
$\prod_{k\in\dN}\spa_k$ defined as follows: 
$x\rF y$ iff $\ans{k:x_k\nE_k y_k}\in\cI$; 

\itla{oe:end}
\label{cp}
{\it the countable power\/} 
$\rF^\iy$ of a \er\ $\rF$ on a space $\spa$ is  
a \er\ on $\spn$ defined as follows: 
$x\rF^\iy y$ iff 
$\ans{[x_k]_{\rF}:k\in\dN}=\ans{[y_k]_{\rF}:k\in\dN},$ 
so that for any $k$ there 
is $l$ with $x_k\rF y_l$ and for any $l$ there 
is $k$ with $x_k\rF y_l$.
\een
Note that operations \ref{cun}, \ref{cdun}, \ref{prod}, 
\ref{cp}, and \ref{fp} with $\cI=\ifi,$ always yield Borel, 
resp., analytic \er s provided given \er s are Borel, 
resp., analytic. 

The operations are not independent. 
In particular, $\bigcap_{k\in\dN}\rF_k$ is Borel reducible 
to $\prod_k\rF_k$ via the map $x\mapsto\ang{x,x,x,...},$ 
while the disjoint union $\bigvee_{k\in\dN}\rF_k$ is 
reducible to $\rav\dN\ti\prod_k\rF_k$ 
via $\ang{k,x}\mapsto \ang{k,x_0,...,x_{k-1},x,x_{k+1},...},$ 
where $x_k\in\spa_k$ are fixed once and for all.
The product $\prod_k\rF_k$ itself is expressible in terms 
of the Fubini product modulo $\ifi.$ 
Indeed, let $f:\dN\onto\dN$ be any map such that $f\obr(n)$ 
is infinite for any $n.$ 
Put $\rE_k=\rF_{f(k)}.$ 
For any $x=\ang{x_0,x_1,x_2,...}\in\prod_k\spa_k$ 
(where $\spa_k$ is the domain of $\rF_k$) let 
$\vt(x)=\ang{y_0,y_1,y_2,...},$ where $y_k=x_{f(k)}.$ 
Then $\vt$ is a Borel reduction of $\prod_k\rF_k$ to 
$\fps{\ifi}{\rE_k}{k}.$ 
Yet the Fubini product and countable power are 
surely not reducible to each other, and we know little 
on the countable union in \ref{cun}. 

It follows that Class 4 of \er s, 
mentioned in our theorems \ref{m1} and \ref{m2}, 
is the least class of \er s which contains equalities 
$\rav\spa$ on Polish spaces $\spa$ and is closed under 
the operations \ref{oe:beg} --- \ref{oe:end} 
(with $\cI=\ifi$ in \ref{fp}), and  
all \er s in Class~4 are Borel \er s on Polish spaces.

There are many interesting \er s in Class~4,  
for instance, the sequence of \er s $\rT_\al\zT\al<\omi,$ 
of H.~Friedman~\cite{frid}  
which begins with 
$\rT_0=\rav\dN,$ the equality relation on $\dN,$ 
then $\rT_{\al+1}={\rT_\al}^\iy$ for any $\al,$ and 
$\rT_\la=\bigvee_{\al<\la}\rT_\ga$ for limit ordinals 
$\la.$ 
Thus, $\dom\rT_1=\dnn$ and 
$x\rT_1 y$ iff $\ran x=\ran y,$ for $x,\,y\in\dnn.$
The map $\vt(x)=$ the characteristic function of 
$\ran x$ witnesses that $\rT_1\reb{\rav\dn}.$ 
To show the converse, define, 
for any $a\in\dn,$ $\ba(a)$ to be the only
increasing bijection $\dN\onto |a|=\ans{k:a(k)=1}$ 
whenewer $|a|$ is infinite, while if 
$|a|=\ans{k_0,...,k_n}$ then put 
$\ba(a)(i)=k_i$ for $i<n$ and $\ba(a)(i)=k_n$ for 
$i\ge n.$ 
The function $\ba$ witnesses ${\rav\dn}\reb{\rT_1},$ 
hence, ${\rT_1}\eqb{\rav\dn}.$  
It easily follows that $\rT_2\eqb {{\rav\dn}}^\iy;$ 
the right-hand side is often taken as the definition 
of the \er\ $\rtd,$ which by this reason is usually 
called ``the equality of countable sets of reals''.

\punk{Ideals}
\label{new:i}

{\it An ideal\/} on a set $A$ is any set 
$\pu\ne\cI\sq\pws A,$ closed under $\cup$ and 
satisfying ${{x\in\cI}\land {y\sq x}}\imp{y\in\cI}.$ 
In this case, $\rei$ is a \er\ on $\pws A$ 
defined as follows: $X\rei Y$ iff ${X\sd Y}\in\cI.$ 
Note that $\rei$ is Borel provided the ideal $\cI$ is such.
Many important \er s appear in the form $\rei,$ among them 
\dm
\bay{rclccl}
\Eo &=& 
\rE_\ifi&,&\text{ where }& 
\ifi=\ans{x\sq\dN:x\,\hbox{ is finite}}\,;\\[\dxii]

\Ei &=& 
\rE_{\Ii}&,&\text{ where }&
\Ii=\fio=
\ans{x\sq\dN^2:\ans{k:\seq xk\ne\pu}\in\ifi}\,;\\[\dxii]

\Et &=& 
\rE_{\It}&,&\text{ where }&
\It=\ofi=\ans{x\sq\dN^2:\kaz k\:(\seq xk\in\ifi)}\,;
\eay
\dm
all of them belong to Class~4~\footnote 
{\ To show that $\Eo$ belongs to Class~4 
let, for any $k,$ $\rF_k$ be the equality on a 
\dd2element set in \ref{fp}. 
To see that $\Et$ belongs to Class~4 
take each $\rF_k$ to be $\Eo$ in \ref{prod}.}.
Ideals of the form 
$\ifi_W=\ans{x\in\pn:x\cap W\in\ifi},$ where 
$W\sq\dN$ is infinite and coinfinite, 
called {\it trivial variations of\/} $\ifi,$ 
also produce \er s in Class~4.

We write $\cI\reb\cJ\zT\cI\eqb\cJ,$ \etc. if resp.\ 
$\rei\reb\rej\zT\,\rei\eqb\rej,$ \etc.

The {\it Fubini product\/} $\fps{\cI}{\cJ_k}{k\in\dN}$ 
of ideals $\cJ_k$ on sets $B_k,$ over an ideal $\cI$ on 
$\dN,$ is an ideal of all sets 
$y\sq B=\ans{\ang{k,b}:k\in \dN\land b\in B_k}$ 
such that the set 
$\ans{k:(y)_k\nin \cJ_k}$ belongs to $\cI,$ where 
$(y)_k=\ans{b:\ang{k,b}\in Y},$ 
a {\it cross-section\/} of $Y.$ 
(Compare with the Fubini product of \er s~!)
In particular, if $\cI,\,\cJ$ are ideals on resp.\ 
$\dN,\,B,$ then $\fpd\cI\cJ=\fps\cI{\cJ_k}{k\in\dN},$ 
where each $\cJ_k=\cJ$ for all $k\in \dN.$   
Thus, $\fpd\cI\cJ$ is the ideal of all sets 
$y\sq \dN\ti B$ such that  
$\ans{k:(y)_k\nin \cJ}\in\cI$.\vom 

{\bfsl P-ideals.\/} 
An ideal $\cI$ on $\dN$ is called {\it a P-ideal\/} if for any 
sequence of sets $x_n\in\cI$ there is a set $x\in\cI$ such 
that $x_n\dif x\in\ifi$ for all $n.$ 
For instance, $\ifi$ and $\It$ (but not $\Ii$!) are P-ideals. 

This class admits different characterizations. 
A {\it submeasure\/} on a set $A$ is any map 
$\vpi:\cP(A)\to[0,\piy],$ satisfying $\vpi(\pu)=0,$ 
$\vpi(\ans a)<\piy$ for all $a,$
and $\vpi(x)\le \vpi(x\cup y)\le \vpi(x)+\vpi(y).$ 
A submeasure $\vpi$ on $\dN$ is 
{\it lover semicontinuous\/}, 
or \lsc\ for brevity, if we have 
$\vpi(x)=\tsup_n\vpi(x\cap\ir0n)$ for all $x\in\pn.$  
Solecki~\cite{sol'}) proved that 
Borel P-ideals are exactly those of the form 
$\Exh_\vpi=\ans{x\in\pn:\vpy(x)=0},$ where $\vpi$ is a 
\lsc\ submeasure on $\dN$ and 
$\vpy(x)=\tinf_n(x\cap\iry n),$ 
and Borel P-ideals are the same as {\it polishable\/} 
ideals, \ie, those which admit a Polish group topology 
with $\sd$ as the group operation. 

Kechris~\cite{rig} proved that the shift action of any Borel 
P-ideal $\cI,$ except for $\ifi\zT\It,$ and trivial 
variations of $\ifi,$ is generically turbulent, 
hence, the corresponding \er\ $\rei$ belongs to Class~1.\vom 

{\bfsl Fr\'echet family.\/} 
This is the least family $\frt$ of ideals containing 
$\ifi$ and closed under the Fubini products 
$\fps{\ifi}{\cI_n}{n\in\dN}.$ 
For instance, the {\it iterated Fr\'echet ideals\/} $\cJ_\al,$ 
defined by induction on $\al<\omi$ so that $\cJ_0=\ifi,$   
$\cJ_{\al+1}=\fpd\ifi{\cJ_\al}$ for all $\al,$ 
and $\cJ_\la=\fps{\ifi_\la}{\cJ_\al}{\al<\la}$ for any 
limit $\la,$ where $\ifi_\la$ is the ideal of all 
finite subsets of $\la,$ belong to $\frt.$  
(A modification of this construction in \cite{jkl}  
involves a cofinal \dd\om sequece fixed in each limit $\la.$)

By definition, if $\cI\in\frt$ then $\rei$ is a \er\ of 
Class~4.

Let $\otp X$ be the order type of $X\sq\Ord.$ 
For any $\ga,\,\al<\omi,$ the set 
\dm
\iz\al\ga \,=\, \ans{A\sq \al:\otp A<\om^\ga} 
\quad\hbox{(non-trivial only if $\al\ge\om^\ga$)}\,.
\dm
is an ideal (because ordinals of the form $\om^\ga$ 
are not sums of a pair of smaller ordinals); 
those ideals, especially in the 
case when $\al=\om^\ga,$ are called {\it indecomposable\/}. 
We don't know whether each $\iz\al\ga$ is really 
(isomorphic to) an ideal in $\frt,$ yet it can be shown 
that any $\iz\al\ga$ is Borel reducible to an ideal in $\frt.$ 
Similarly, 
\dm
\iw\al\ga \,=\, \ans{A\sq \ga:\rkb A<\om^\ga}  
\quad\hbox{(the {\it Weiss ideals\/}, 
non-trivial only if $\al\ge\omm\ga$)},
\dm
where $\rkb X$ is the Cantor-Bendixson rank of $X\sq\Ord,$ 
see Farah~\cite[\S~1.14]{aq}, 
are Borel reducible to ideals in $\frt$.

\parf{Proof of Theorem~\ref{m1}}
\label{thm1}

This section proves Theorem~\ref{m1}. 
Our method of demonstration of the non-reducibility employs 
the following auxiliary notions. 
Let $\rE,\,\rF$ be \er s on Polish spaces, 
resp., $\dX,\,\dY.$ 
A map $\vt:\dX\to\dY$ is called
\bit
\item\msur
\ef{\it invariant\/} \ if \ 
${x\rE y}\imp {\vt(x)\rF\vt(y)}$ for all 
$x,\,y\in\dX$; 

\item
{\it generically\/}~\footnote
{\ In this research direction, ``generically'', or, in our 
abbreviation, ``\genp'' (property) is understood so  
that the property holds on a comeager domain.} 
(or {\it\genp\/}, for brevity) 
\ef{\it invariant\/} if ${x\rE y}\imp {\vt(x)\rF\vt(y)}$ 
for all $x,\,y$ in a comeager set $X\sq\dX$;

\item 
\hspace{0pt}{\it\gen \ddf constant\/} \ if \ 
${\vt(x)\rF\vt(y)}$ holds for all $x,\,y$ 
in a comeager set~$X\sq\dX.$
\eit 
Finally, $\rE$ is {\it generically\/} 
\ddf{\it ergodic\/} (Hjorth~\cite[3.1]h) 
if every \bm\ \ef invariant map is \gen \ddf constant. 

\bpro
\label{t2nonr}
${\rm(i)}$ 
If\/ $\rE$ is \gen \ddf ergodic and 
does not have a comeager equivalence class  
then\/ $\rE$ is not reducible to\/ $\rF$ by a 
\bm\ map.

${\rm(ii)}$ 
If\/ $\rE$ is \gen \ddf ergodic then even every \bm\ 
\gen\ef invariant map is \gen \ddf constant.\qed
\epro
 
Accordingly, our plan will be to show that any orbit 
\er\ induced by a Polish turbulent action is \gen \ddf ergodic 
for any \er\ $\rF$ in Class~4.

\punk{Local orbits and turbulence}
\label{tu:def}

An {\it action\/} of a group $\dG$ on $\dX$ is 
any map $a:{\dG\ti\dX}\to\dX,$ usually written as 
$a(g,x)=g\app x,$ such that 
$1)\msur$ $e\app x=x,$ and 
$2)\msur$ $g\app(h\app x)=(gh)\app x.$ 
In this case, $\stk\dX a,$ as well as $\dX$ itself, 
is called a {\it\dd\dG space\/}.  
A continuous action of a Polish group~\footnote
{\ That is, a topological group whose underlying set is a 
Polish space and the group operation and the inverse map 
are continuous.}  
$\dG$ on a Polish space $\dX$ is called {\it Polish action\/}  
and $\dX$ itself is called a {\it Polish\/ \dd\dG space\/}. 

Any action $a$ of $\dG$ on $\dX$ induces 
{\it the orbit \er\/},  
$\aer a\dX=\ergx,$ defined on $\dX$ so that $x\ergx y$ 
iff there is $a\in\dG$ with $y=a\app x,$ whose equivalence 
classes 
\dm
[x]_\dG=[x]_{\ergx}=\ans{y:\sus g\in\dG\;(g\app x=y)}\,.
\dm
are \dd\dG orbits. 
Induced \er s of Polish actions  
are analytic (as sets of pairs), sometimes even 
Borel \cite[Chapter 7]{beke}.

Suppose that a group $\dG$ acts on a space $\dX.$ 
If $G\sq\dG$ and $X\sq\dX$ then we define 
\dm
\rr XG=\ans{\ang{x,y}\in X^2:\sus g\in G\:(x=g\ap y)}
\dm 
and let $\sym XG$ denote the \er-hull of $\rr XG,$ \ie, 
the \dd\sq least \er\ on $X$ such that 
${x\rr XG y}\imp {x\sym XG y}.$ 
In particular ${\sym\dX\dG}={\ergx},$ but generally 
we have ${\sym XG}\sneq{\ergx\res X}.$  
Finally, define 
$\lo XGx=\ek x{\sym XG}=\ans{y\in X:x\sym XG y}$ 
for $x\in X$ -- the {\it local orbit\/} of $x.$ 
In particular, $\ek x\dG=\ek x{\ergx}=\lo\dX\dG x,$ 
the full \dd\dG orbit of $x\in\dX$.

\bdf[{{\rm This version is taken from Kechris~\cite[\S~8]{apg}}}]
\label{df:durb}
Suppose that $\dX$ is a Polish space and $\dG$ is 
a Polish group acting on $\dX$ continuously. 
\ben
\tenu{(t\arabic{enumi})}
\itla{T1}
A point $x\in\dX$ is {\it turbulent\/} if for any open
non-empty set $X\sq \dX$ containing $x$ and any nbhd  
$G\sq\dG$ (not necessarily a subgroup) of  
$\ong,$ the local orbit $\lo XGx$ is somewhere 
dense (\ie, not a nowhere dense set) in $\dX$.

\itla{T2} 
An orbit $\ek x\dG$ is {\it turbulent\/} if $x$ is such 
(then all $y\in\ek x\dG$ are turbulent).

\itla{T3}
The action (of $\dG$ on $\dX$) is {\it\gen turbulent\/} 
and $\dX$ is a {\it \gen turbulent\/} Polish \dd\dG space, 
if all orbits are dense and meager. 
The action is {\it generically}, 
or {\it\gen turbulent\/} and $\dX$ is a 
{\it \gen turbulent\/} Polish \dd\dG space, if 
the union of all dense, turbulent, and meager orbits 
$[x]_\dG$ is comeager.\qed
\een
\eDf

\er s induced by \gen turbulent Polish actions are what is 
called {\it Class 1\/} in Theorem~\ref{m1}. 
To prove Theorem~\ref{m1}, we are going to show that 
all \er s in Class~1 are \gen \ddf ergodic, for any $\rF$ 
in Class~4. 
The method of proof will be by induction on the 
construction of \er s in Class~4. 
It is a slight inconvenience that we have to consider a 
somewhat stronger property through the induction scheme. 

Suppose that $\rF$ is an \er\ on a Polish space. 
\bit
\item
An action of $\dG$ on $\dX$ 
is {\it hereditarily generically\/} 
({\it\hgp,} for brevity) \ddf ergodic   
if \er\ $\sym XG$ is generically \ddf ergodic  
whenever $X\sq \dX$ is a non-empty open set,  
$G\sq\dG$ is a non-empty open set containing $\ong,$ 
and local orbits $\lo{X}Gx$ are 
dense in $X$ for comeager (in $X$) many $x\in X$. 
\eit
This obviously implies \gen \ddf ergodicity  
provided the action is \gen turbulent. 
 
\bte
\label{Turb}
Suppose that\/ $\dG$ is a Polish group, $\dX$ is a 
\gen turbulent Polish\/ \dd\dG space. 
Then\/ $\ergx$ is\/ \hg\ddf ergodic, hence\/ 
{\rm(by Proposition~\ref{t2nonr})}, \poq{not} 
reducible to any \er\ $\rF$ of Class~4 by a BM map.
\ete

\punk{Preliminaries to the proof of Theorem~\ref{Turb}}
\label{Tprel}

We begin with two rather simple technical facts related 
to turbulence. 

\ble
\label{deno}
In the assumptions of the theorem, suppose that 
$\pu\ne X\sq\dX$ is an open set, $G\sq\dG$ is a nbhd 
of\/ $\ong,$ and\/ $\lo XGx$ is dense in\/ $X$ for\/ 
\dd Xco\-meager many\/ $x\in X.$ 
Let\/ $U,\,U'\sq X$ be non-empty open and\/ $D\sq X$ 
be comeager in\/ $X.$ 
Then there exist points\/ $x\in D\cap U$ and\/ 
$x'\in D\cap U'$ with\/ $x\sym XG x'$.
\ele
\bpf
Under our assumptions there exist points $x_0\in U$ and 
$x'_0\in U'$ with $x_0\sym XG x'_0,$ \ie, there are elements 
$g_1,...,g_n\in G$ with 
$x'_0=g_n\app g_{n-1}\app...\app g_1\app x_0$ and  
$g_k\app ...\app g_1\app x_0\in X$ for all $k\le n.$ 
Since the action is continuous, there is a nbhd $U_0\sq U$ of 
$x_0$ such that $g_k\app ...\app g_1\app x\in X$ 
and $g_n\app g_{n-1}\app...\app g_1\app x\in U_2$ for \poq{all}  
$x\in U_0.$ 
Since $D$ is comeager, easily there is $x\in U_0\cap D$ such 
that  $x'=g_n\app g_{n-1}\app...\app g_1\app x\in U'\cap D$. 
\epf

\ble
\label{kl}
In the assumptions of the theorem, 
for any open non-empty\/ $U\sq\dX$ and\/ 
$G\sq\dG$ with\/ $\ong\in G$ there is an open set\/ 
$\pu\ne U'\sq U$ such that the local orbit\/ $\lo{U'}Gx$ 
is dense in $U'$ for\/ \dd{U'}comeager many\/ $x\in U'$.
\ele
\bpf
Let $\incl X$ be the interior of the closure of $X.$ 
If $x\in U$ and $\lo UGx$ is somewhere dense (in $U$) 
then the set $U_x=U\cap\incl{\lo UGx}\sq U$ is open and\/ 
\dd{\sym UG}invariant 
(an observation made, \eg, in \cite[proof of 8.4]{apg}), 
moreover, ${\lo UGx}\sq U_x,$ hence, ${\lo UGx}={\lo{U_x}Gx}.$ 
It follows from the invariance that the sets $U_x$ are 
pairwise disjoint, and it follows from the turbulence that 
the union of them is dense in $U.$  
Take any non-empty $U_x$ as $U'.$ 
\epf
 
Our proof of Theorem~\ref{Turb} goes on by induction on the 
construction of \er s in Class~4 with the help of 
the operations mentioned in Section~\ref{new:er}, 
in several following subsections. 
We begin with the base of the induction: prove that, under the 
assumptions of the theorem, $\ergx$ is \hg\dd{\rav\dN}ergodic. 
Suppose that $X\sq\dX$ and $G\sq\dG$ are non-empty open 
sets, $\ong\in G,$ and $\lo{X}Gx$ is dense in $X$ for 
\dd Xcomeager many points $x\in X.$ 
Prove that $\sym XG$ is 
generically \dd{\rav\dN}ergodic. 

Consider a \bm\  
\gen\dd{\sym XG,\rav\dN}invariant map $\vt:\dX\to\dN.$ 
Suppose, on the contrary, that $\vt$ is not 
\gen \dd{\rav\dN}constant. 
Then there exist open non-empty sets $U_1,\,U_2\sq X,$ 
numbers $\ell_1\ne \ell_2,$ and a comeager set $D\sq X$ 
such that $\vt(x)=\ell_1$ for all $x\in D\cap U_1$ and 
$\vt(x)=\ell_2$ for all $x\in D\cap U_2.$ 
Lemma~\ref{deno} yields a pair of points 
$x_1\in U_1\cap D$ and 
$x_2\in U_2\cap D$ with $x_1\sym XG x_2,$ contradiction. 

\punk{Inductive step of the countable power}
\label{tut:iy}

Consider a \gen turbulent Polish \dd\dG space $\dX$  
and a Borel \er\ $\rF$ on a Polish space $\dY.$ 
Assume that the action of $\dG$ on $\dX$ is \hg\ddf ergodic, 
and prove that the action is \hg\dd\rfy ergodic. 
Fix a non-empty open set $X_0\sq\dX$ and a nbhd $G_0$ of $\ong$ 
in $\dG,$ such that $\lo{X_0}{G_0}x$ is dense in $X_0$ for 
\dd{X_0}comeager many $x\in X_0.$  
Prove that any given \dd{\sym{X_0}{G_0},\rfy}invariant \bm\ 
function $\vt:X_0\to\dY^\dN$ is \gen \dd\rfy constant.
By definition, we have 
\ben
\tenu{(\arabic{enumi})}
\itla{X1} 
for $x,\,x'\in X_0:$  
$x\sym{X_0}{G_0} x'\limp
\kaz k\;\sus l\;\skl\vt_k(x)\rF\vt_l(x')\skp,$ 
where $\vt_k(x)=\vt(x)(k)$. 
\een
Note that $\vt$ is continuous on a dense $\Gd$ set $D\sq X_0$. 

\ble
\label{L1}
For each\/ $k$ and open\/ $\pu\ne U\sq X_0$ 
there is an open set\/ $\pu\ne W\sq U$ such that\/ 
$\vt_k$ is \gen \ddf constanta on\/ $W$.
\ele
\bpf
A simple category argument beginning with \ref{X1} yields a 
number $l$ and open non-empty sets $W\sq U$ and $Q\sq G_0,$ 
and a dense in $W\ti Q$ set $P\sq W\ti Q$ of class $\Gd$ 
such that ${\vt_k(x)}\rF{\vt_l(g\ap x)}$ holds for all pairs 
$\ang{x,g}\in P.$ 
We can assume that ${\ang{x,g}\in P}\imp {x\in D}.$ 
Since $Q$ is open, there is an element 
$g_0\in Q$ and a nbhd $G\sq G_0$ of $\ong$ 
with $G\obr=G$ such that $g_0\,G\sq Q$.

The continuation of the proof involves forcing~\footnote
{\  
Some degree of the reader's acquaintance with forcing is 
assumed. 
The lemma could have been proved by purely topological 
arguments, 
yet then the reasoning then would not be so transparent.}. 

Let us fix a countable transitive model 
$\mm$ of $\zhc,$ 
\ie, $\ZFC$ minus the Power Set axiom but plus the axiom: 
``every set is hereditarily countable''. 
We can assume that
$\dX$ is coded in $\mm$ in the sense that there is a  
set $D_\dX\in\mm$ which is a dense (countable) subset of 
$\dX,$ and $d_\dX\res D_\dX$ 
(the distance function of $\dX$ restricted to $D_\dX$) 
also belongs to $\mm.$ 
Further, 
$\dG\zT\dY,$ the action of $\dG$ on $\dX,$ sets $G,\,D,\,P$ 
and the map $\vt\res{D}$ 
are also assumed to be coded in $\mm$ in a similar sense. 

Below, let $\px$ be the Cohen forcing for $\dX,$ 
which consists of rational balls with centers in  
a fixed dense countable subset of $\dX,$ and let $\pg$ 
be the Cohen forcing for $\dG$ defined similarly 
As usual, $U\sq V$ means that $U$ is a stronger 
forcing condition.  
In these assumptions, the notions of Cohen generic, over 
$\mm,$ elements of $\dX$ and $\dG$ makes sense, and the 
set of all Cohen generic, over $\mm$ points of $\dX$ is a 
dense $\Gd$ subset of $\dX$ included in $D$. 

\bcl[{{\rm The key point of the turbulence}}]
\label{41*}
If\/ $x,\,x'\in W$ are\/ \dd\px generic over\/ $\mm$  
and\/ $x\sym{W}{G}x'$ then\/ $\vt_k(x)\rF\vt_k(x')$.
\ecl
\bpf 
We argue by induction on the number $n(x,x')$ 
equal to the least $n$ 
such that there exist $g_1,...,g_n\in G$ satisfying
\ben
\tenu{(\arabic{enumi})}
\aci1
\itla{X2}
$x'=g_n\app g_{n-1}\app...\app g_1\app x,$ \ and \
$g_k\app...\app g_1\app x\in {W}$ \ for all \ $k\le n$.
\een
Suppose that $n(x,x')=1,$ thus, $x=h\ap x'$ for some 
$h\in G\cap\mm[x,x']$~\footnote
{\label{mxy}\ 
Here $\mm[x,x']$ is defined as any (countable transitive) 
model 
of $\zhc$ containing $x,\,x',$ and all sets in $\mm,$ rather 
than a generic extension of $\mm.$ 
The model $\mm[x,x']$ can contain more ordinals than 
$\mm,$ but this is not essential here.}
Take any \dd\pg generic, over $\mm[x,x'],$ element $g\in Q,$ 
close enough to $g_0$ for $g'=gh\obr$ to belong to $Q.$ 
Then $\ang{x,g}$ is \dd{\px\ti\pg}generic over $\mm$ 
by the product forcing theorem, thus, 
$\ang{x,g}\in P$ 
(because $P$ is a dense $\Gd$ coded in $\mm$) 
and $\vt_k(x)\rF\vt_l(g\ap x)$ by the choice of $P.$ 
Moreover, $g'$ also is \dd{\pg}generic over $\mm[x'],$ 
so that $\vt_k(x')\rF\vt_l(g'\ap x')$ by the same argument. 
Yet we have $g'\ap x'=gh\obr\ap(h\ap x)={g\ap x}$. 

As for the inductive step, suppose that \ref{X2} holds 
for some $n\ge 2.$ 
Take a \dd\pg generic, over $\mm[x],$ element $g'_1\in G$ 
close enough to $g_1$ for $g'_2= g_2\,g_1\,{g'_1}\obr$ to 
belong to $G$ and for $x^\ast=g'_1\app x$ to belong to $W.$ 
Note that $x^\ast$ is \dd\px generic over $\mm$ 
(product forcing) 
and $n(x^\ast,x')\le n-1$ because 
$g'_2\app x^\ast=g_2\app g_1\app x$.\vom
 
\epF{Claim~\ref{41*}}

To summarize, we have shown that $\vt_k$ is  
\gen\dd{\sym{W}{G},\rF}invariant on $W,$ 
\ie, invariant on a comeager subset of $W.$  
We can also assume that the orbit $\lo{W}Gx$ 
is dense in $W$ for \dd{W}comeager many points $x\in W,$ 
by Lemma~\ref{kl}. 
Then, by the \hg\ddf ergodicity, $\vt_k$ is 
\gen \ddf constant on $W,$ as required.

\epF{Lemma~\ref{L1}}

According to the lemma, there exist: an \dd{X_0}comeager set 
$Z\sq X_0,$ and a countable set $Y=\ans{y_j:j\in\dN}\sq \dY$ 
such that, for any $k$ and for any $x\in Z$ there is $j$ with  
$\vt_k(x)\rF y_j.$ 
Let $\eta(x)=\bigcup_{k\in\dN}\ans{j:\vt_k(x)\rF y_j}.$ 
Then, for any pair $x,\,x'\in Z,$ 
$\vt(x)\rfy \vt(x')$ iff $\eta(x)=\eta(x'),$ so that, by 
the invariance of $\vt,$ we have: 
\ben
\tenu{(\arabic{enumi})}
\aci2
\itla{X3}
${x\sym{X_0}{G_0}x'}\,\limp \,{\eta(x)=\eta(x')}$ \ 
for all $x,\,x'\in Z$.
\een
It remains to show that $\eta$ is a constant on a comeager 
subset of $Z$.  

Suppose, on the contrary, that there exist two non-empty 
open sets $U_1,U_2\sq X_0,$ a number $j\in\dN,$ and a 
comeager set $Z'\sq Z$ such that 
$j\in\eta(x_1)$ and $j\nin\eta(x_2)$ for all 
$x_1\in Z'\cap U_1$ and $x_2\in Z'\cap U_2.$ 
Lemma~\ref{deno} yields a contradiction to \ref{X3}
as in the end of Section~\ref{Tprel}.

\punk{Inductive step of the Fubini product}
\label{tut:Fu}

Suppose that $\dX$ is a \gen turbulent Polish \dd\dG space. 
Prove that the action of $\dG$ on $\dX$ is 
\hg\dd\rF ergodic, where $\rF=\fps{\ifi}{\rF_k}{k},$
$\rF_k$ is a Borel \er\ on a Polish space $\dY_k,$ and 
the action is \hg\dd{\rF_k}ergodic for any $k.$ 
Fix an open set $\pu\ne X_0\sq\dX$ and a nbhd $G_0$ of 
$\ong$ in $\dG,$ such that \dd{X_0}comeager many orbits 
$\lo{X_0}{G_0}x$ with $x\in X_0$ are dense in $X_0.$   
Prove that any \dd{\sym{X_0}{G_0},\rF}invariant \bm\ 
function $\vt:U_0\to\dY$ is \gen \ddf constant on $X_0.$ 
By definition  
\ben
\tenu{(\arabic{enumi})}
\aci3
\itla{Y1}
for $x,\,y\in X_0:$ 
$x\sym{X_0}{G_0} y\limp 
\sus k_0\;\kaz k\ge k_0\;\skl\vt_k(x)\rF_k\vt_k(y)\skp$, 
\een 
where $\vt_k(x)=\vt(x)(k).$ 
Note that $\vt$ is continuous on a dense $\Gd$ set $D\sq X_0$.

\ble
\label{L2}
For any open set\/ $\pu\ne U\sq X_0$ 
there exist a number\/ $k_0$ and open\/ $\pu\ne W\sq U$   
such that\/ $\vt_k$ is \gen \ddf constant on\/ $W$ 
for all\/ $k\ge k_0$.
\ele
\bpf
Applying \ref{Y1}, we can easily find a number $k_0,$  
open non-empty sets $W\sq U$ and $Q\sq G_0,$ and a dense  
in $W\ti Q$ set $P\sq W\ti Q$ of class $\Gd,$ 
such that ${\vt_k(x)\rF\vt_k(g\ap x)}$ holds for all 
$k\ge k_0$ and all pairs $\ang{x,g}\in P.$ 
We can assume that ${\ang{x,g}\in P}\imp {x\in D}.$ 
Since $Q$ is open, there exist an element 
$g_0\in Q$ and a nbhd $G\sq G_0$ of $\ong$ 
with $G\obr=G,$ such that $g_0\,G\sq Q$. 

Let a model $\mm$ be as in the proof of  
Lemma~\ref{L1}. 
Similarly to Claim~\ref{41*}, we can prove that  
if points $x,\,x'\in W$ are  
\dd\px generic over $\mm,$ $k\ge k_0,$ and  
$x\sym{W}{G}x',$ then $\vt_k(x)\rF_k\vt_k(x'),$ 
in other words, each function $\vt_k$ with $k\ge k_0$ 
is \gen\dd{\sym{W}{G},\rF_k}invariant on $W.$ 
We can assume, by Lemma~\ref{kl}, that  
\dd{W}comeager many orbits $\lo{W}Gx$ 
are dense in $W.$  
Now, by the \hg\dd{\rF_k}ergodicity, any $\vt_k$ with 
$k\ge k_0$ is \gen \dd{\rF_k}constant on such a set $W,$ 
as required.
\epf

It is clear that if $W$ is chosen as in the lemma then  
$\vt$ itself is \gen \dd{\rF}constant on $W.$ 
It remains to show that these constants are  
\ddf equivalent to each other. 
Suppose, on the contrary, that there exist two non-empty 
open sets $W_1,\,W_2\sq X_0$ and a pair of points $y\nF y'$ 
in $\dY$ such that $\vt(x)\rF y$ and $\vt(x')\rF y'$ for 
comeager many $x\in W_1$ and $x'\in W_2.$ 
Contradiction follows as in the end of 
Section~\ref{tut:iy}.

\punk{Other inductive steps}
\label{tut:oth}

We carry out the induction steps related to operations 
\ref{cun}, \ref{cdun}, \ref{prod} of \ref{new:er}.\vom

{\sl Countable union\/}. 
Suppose that $\rF_1,\,\rF_2,\,\rF_3,\,...$ are Borel \er s   
on a Polish space\/ $\dY,$ and $\rF=\bigcup_k{\rF_k}$ is 
still a \er, and the Polish and \gen turbulent action of 
$\dG$ on $\dX$ is \hg\dd{\rF_k}ergodic for any $k,$ 
and prove that it remains \hg\ddf ergodic. 

Consider a non-empty open set $X_0\sq\dX$ and a nbhd $G_0$ 
of $\ong$ in $\dG,$ such that \dd{X_0}comeager many orbits 
$\lo{X_0}{G_0}x$ with $x\in X_0$ are dense in $X_0.$   
Consider a \dd{\sym{X_0}{G_0},\rF}invariant \bm\ function 
$\vt:X_0\to\dY,$ continuous on a dense $\Gd$ set $D\sq X_0.$  
It follows from the invariance that for any open set 
$\pu\ne U\sq X_0$ there exist: a number $k$ and open 
non-empty sets $Q\sq U$ and $Q\sq G_0$  
such that $\vt(x)\rF_k\vt(g\ap x)$ holds
for any \dd{\px\ti\pg}generic, over $\mm,$    
pair $\ang{x,g}\in W\ti Q.$  
We can find, as above, $g_0\in Q\cap\mm$ 
and a nbhd $G\sq G_0$ of $\ong$ such that $g_0\,G\sq Q.$
Similarly to Claim~\ref{41*}, we have 
$\vt(x)\rF_k\vt(x')$ for any pair of  
\dd\px generic, over $\mm,$ elements $x,\,x'\in W,$ 
satisfying $x\sym{W}{G}x'.$
It follows, by the ergodicity, that 
$\vt$ is \gen\dd{\rF_k}constant, hence, \ddf constant,  
on $W.$ 
That these \ddf constants are \ddf equivalent to each other, 
can be demonstrated exactly as in the end of 
Section~\ref{tut:iy}. 
The operation of countable intersection is considered 
similarly.\vom

{\sl Countable product\/}. 
It is shown in \ref{new:er} that this operation is reducible  
to the Fubini product, yet there is a simple independent  
argument.
If $\rF_k$ be \er s on spaces $\dY_k$ 
then $\rF=\prod_k{\rF_k}$ is a \er\ on the space 
$\dY=\prod_k\dY_k.$ 
For any map $\vt:\dX\to\dY,$ to be \dd{\rE,\rF}invariant 
(where $\rE$ is an arbitrary \er\ on $\dX$) 
it is necessary and sufficient that every 
co-ordinate map $\vt_k(x)=\vt(x)(k)$ 
is \dd{\rE,\rF_k}invariant, which immediately 
yields the result required.
\vom

{\sl Disjoint union\/}. 
It is shown in \ref{new:er} that this operation is reducible  
to the product.\vtm

\qeDD
{Theorems~\ref{Turb} and \ref{m1}}

\parf{Applications}
\label{appl}

This section contains two applications of Theorem~\ref{Turb}. 
One of them is Theorem~\ref{m2}. 
The other one shows how Hjort's theorem mentioned in the 
Introduction (that ``turbulent'' \er s are not Borel 
reducible to Polish actions of $\isg,$ the group of all 
permutations of $\dN$) can be derived 
from Theorem~\ref{Turb} by rather simple arguments.

\punk{Proof of Theorem~\ref{m2}}
\label{thm2}

Let us fix a non-trivial, as in Theorem~\ref{m2},  
Borel P-ideal $\cZ\sq\pn.$  
By a theorem of Solecki (see Section~\ref{new:i}) 
there exists a \lsc\ submeasure $\vpi$ on $\dN$ such that 
$\cZ=\ans{x\sq\dN:\vpy(x)=0}.$  
Put $r_k=\vpi(\ans k)$.

\ble[{{\rm Kechris~\cite{rig}}}]
\label{p2s}
If\/ $\cZ$ is not equal to\/ $\ifi,$ is not a trivial 
variation of\/ $\ifi,$ and is not isomorphic to\/ $\cI_3=\ofi,$  
then there is a set\/ $W\nin\cZ$ such that\/ 
$\sis{r_k}{k\in W}\to0$.
\ele
\bpf  
Put   
$U_n=\ans{k:r_k\le\frac1n},$ separately 
$U_0=\dN,$ thus, $U_{n+1}\sq U_n$ for all $n.$ 
We claim that $\tinf_{m\in\dN}\vpi(U_m)>0.$ 
Otherwise a set $x\sq\dN$ belongs to $\cZ$ iff 
$x\dif U_n$ is finite for any $n.$ 
If the set  
$N=\ans{n:U_n\dif U_{n+1}\,\text{ is infinite}}$ 
is empty then easily $\cZ=\pn.$ 
If $N\ne\pu$ is finite then $\cZ$ is either $\ifi$ 
(if eventually $U_n=\pu$)
or a trivial variation of $\ifi$ 
(if $U_n$ is non-empty for all $n$).
If finally $N$ is infinite then $\cZ$ 
is isomorphic to $\ofi.$ 
(For instance, if all sets $D_n=U_n\dif U_{n+1}$ are 
infinite then $x\in\cZ$ iff $x\cap D_n$ is finite for 
all $n$.)  
Thus we always have a contradiction to the 
assumptions of the lemma.

It follows that there is $\ve>0$ such that $\vpi(U_m)>\ve$ 
for all $m.$ 
As $\vpi$ is \lsc, we can define an increasing sequence 
of numbers $n_1<n_2<n_3<...$ and for any $l$ a finite set 
$w_l\sq U_{n_l}\dif U_{n_{l+1}}$ with $\vpi(w_l)>\ve.$ 
Then $W=\bigcup_lw_l\nin\cZ$ and obviously 
$\ans{r_k}_{k\in W}\to 0$. 
\epf

Since obviously $\rE_{\cZ\res W}\reb\rez,$ the following 
lemma is sufficient for Theorem~\ref{m2}:

\ble
\label{ket}
If\/ $\cZ\zT\vpi\zT r_k$ are as above, and\/ 
$\sis{r_k}{}\to0,$ then the shift action of\/ $\cZ$ on\/ 
$\pn$ is \gen turbulent. 
\ele
\bpf
$\cZ$ is a Polish group (with the operation $\sd$) in the 
topology $\tau$ induced by the metric $r(x,y)= \vpi(x\sd y).$ 
The action of $\cZ$ by $\sd$ on the space $\pn$ 
(considered in the product topology; $\pn$ is here 
identified with $\dn$) 
by $\sd$ is then continuous. 
It remains to verify the turbulence. 

Let $x\in\pn.$  
The orbit $\ek x\cZ=\cZ\sd x$ is easily dense and meager, 
hence, it suffices to prove that $x$ is a turbulent point 
of the action. 
Consider an open set $X\sq\pn$ containing $x,$ and 
a \dd\tau hbhd $G$ of $\pu$ (the neutral element of $\cZ$); 
we may assume that, for some $k,$ 
$X=\ans{y\in\pn:y\cap [0,k)=u},$ where $u=x\cap[0,k),$ and 
$G=\ans{g\in\cZ:\vpi(g)<\ve}$ for some $\ve>0.$ 
Prove that the local orbit $\lo XGx$ is somewhere 
dense (\ie, not a nowhere dense set) in $X$. 

Let $l\ge k$ be big enough for $r_n<\ve$ for all $n\ge l.$ 
Put $v=x\cap [0,l)$ and prove that $\lo XGx$ is dense in 
$Y=\ans{y\in\pn:y\cap [0,l)=v}.$ 
Consider an open set  
$Z=\ans{z\in Y:z\cap [l,j)=w},$ where $j\ge l\zT w\sq[l,j).$ 
Let $z$ be the only element of $Z$ with 
$z\cap{[j,\piy)}=x\cap{[j,\piy)},$ thus, 
${x\sd z}=\ans{l_1,...,l_m}\sq[l,j).$ 
Each $g_i=\ans{l_i}$ belongs to $G$ by the choice of $l$ 
(indeed, $l_i\ge l$). 
Moreover, easily 
$x_i=g_i\sd g_{i-1}\sd...\sd g_1\sd x=
\ans{l_1,...,l_i}\sd x$ 
belongs to $X$ for any $i=1,...,m,$ and $x_m=z,$ 
thus, $z\in \lo XGx,$ as required.\vtm

\epF{Lemma and Theorem~\ref{m2}}

\punk{Irreducibility to actions of the 
group of all permutations of $\dN$}
\label{groas}
 
Recall that $\isg$   
is the group of all permutations of $\dN,$ \ie, 1--1 maps 
$\dN\onto\dN,$ with the superposition as the group operation. 
A compatible Polish metric on $\isg$ can be defined by 
$D(x,y)=d(x,y)+d(x\obr,y\obr),$ where $d$ is the ordinary 
Polish metric of $\dnn,$ \ie, $d(x,y)=2^{-m-1},$ where 
$m$ is the least number such that $x(m)\ne y(m)$. 

Hjorth 
proved in mid-90s that turbulent \er s are not reducible to 
those induced by Polish actions of $\isg.$ 
The proof (as, \eg, in \cite{h,apg}) is quite 
complicated, in particular, containing references to some 
model theoretic facts and methods like Scott's analysis. 
We decided to include a simplified proof, based on the 
following theorem. 
This will be still a lengthy argument, 
because, to make the exposition friendly to a reader not 
experienced in special topics related to group actions 
and model theory, we outline proofs of some 
auxiliary results involved.

\bte
\label H
Any \er\/ $\rE,$ induced by a Polish action 
and reducible to an orbit \er\ of a Polish action 
of\/ $\isg$ by a \bm\ map, is reducible to one of \er s\/ 
$\rT_\ga$ by a \bm\ map~\footnote
{\ We cannot claim Borel reducibility here, because, as all 
\er s $\rT_\ga$ are easily Borel, any \er\ Borel reducible 
to some $\rT_\ga$ is Borel itself, while 
on the other hand even \er s of the 
form $\ism\cL$ are, generally, non-Borel (but analytic).} 
--- hence, by Theorem~\ref{m1}, 
such an \er\/ $\rE$ cannot be induced by a \gen turbulent 
Polish action.
\ete

\punk{Classifiability by countable structures}
\label{loac}

Isomorphism relations of various classes of countable 
structures are amongst those induced by Polish actions of 
$\isg.$ 
Indeed, suppose that $\cL=\sis{R_i}{i\in I}$ is a countable 
relational language, \ie, $\card I\le\alo$ and each $R_i$ 
is an \dd{m_i}ary relational symbol. 
Put~\footnote
{\ $X_\cL$ is often used to denote $\mox\cL$.}   
$\mox\cL=\prod_{i\in I}\cP(\dN^{m_i}),$ 
the space of \dd\cL{\it structures\/} on $\dN$ 
as the underlying set.
{\it The logic action\/} $\loa\cL$ of $\isg$ on $\mox\cL$ 
is defined as follows: 
if $x=\sis{x_i}{i\in I}\in\mox\cL$ and $g\in\isg$ 
then $y=\loa\cL(g,x)=g\app x=\sis{y_i}{i\in I}\in\mox\cL,$ 
where 
\dm
{\ang{k_1,...,k_{m_i}}\in x_i} \leqv 
{\ang{g(k_1),...,g(k_{m_i})}\in y_i}  
\dm
for all $i\in I$ and $\ang{k_1,...,k_{m_i}}\in\dN^{m_i}.$ 
Then $\stk{\mox\cL}{\loa\cL}$ is a Polish \dd\isg space,   
while \dd{\loa\cL}orbits in $\mox\cL$ are exactly the 
isomorphism classes of \dd\cL structures, which is a reason 
to denote the associated equivalence relation 
$\aer{\loa\cL}{\mox\cL}$ by ${\ism\cL}.$ 
All \er s of the form $\ism\cL$ are analytic, of course. 

Hjorth~\cite[2.38]h defined an \er\ $\rE$ to be 
{\it classifiable by countable structures\/} 
\index{ER!classifiable by countable structures}%
if there is a countable relational language $\cL$ such that 
$\rE\reb{\ism\cL}$.

\bte[{{\rm Becker and Kechris \cite{beke}}}]
\label{siy2cc}
Any \er\ induced by a Polish action of\/ 
$\isg$ is classifiable by countable structures. 
\ete

Thus all \er s induced by Polish actions of $\isg$ 
(in fact also of any closed subgroup of $\isg$) are 
Borel reducible to a very special kind of actions of $\isg.$ 

\bpf[{{\rm by Hjorth~\cite[6.19]h}}]  
Let $\dX$ be a Polish \dd\isg space with basis 
$\sis{U_l}{l\in\dN},$ and let $\cL$ be the language with 
relations $R_{lk}$ of arity $k.$ 
If $x\in \dX$ then define $\vt(x)\in\mox\cL$ by 
stipulation that $\vt(x)\mo R_{lk}(s_0,...,s_{k-1})$ 
if and only if $s_i\ne s_j$ whenever $i<j<k,$ and   
$g\obr\app x\in U_l$ whenever  
$g\in\isg$ satisfies $\ang{s_0,...,s_{k-1}}\su g.$  
Then $\vt$ reduces $\aer\isg\dX$ to $\ism\cL.$ 
\epf

\punk{Reduction to countable graphs}
\label{gra}

It could be expected that more complicated languages  
$\cL$ produce more complicated \er\ $\ism\cL.$ 
However this is not the case: 
it turns out that a single binary relation can code 
structures of any countable language.  
Let $\cG$ be the language of (oriented binary) graphs, \ie, 
$\cG$ contains a single binary predicate, say $R(\cdot,\cdot)$.

\bte
\label{grafs}
If\/ $\cL$ is a countable relational language then\/ 
${\ism\cL}\reb{\ism\cG}$. 
\ete  

Becker and Kechris \cite[6.1.4]{beke} outline a proof based 
on coding in terms of lattices, unlike the following argument, 
yet it may in fact involve the same idea.

\bpf
Let $\hfn$ be the set of all hereditarily finite sets over the 
set $\dN$ considered as the set of atoms, and $\ve$ be the 
associated ``membership'' 
(no $n\in\dN$ has  
\dd\ve elements, $\ans{0,1}$ is different from $2,$ \etc\/.). 
Let $\ihf$ be the $\hfn$ version of $\ism\cG,$ \ie, if 
$P,\,Q\sq\hfn^2$ then $P\ihf Q$ means that there is a 
bijection $b$ of $\hfn$ on itself such that 
$Q=b\app P=\ans{\ang{b(s),b(t)}:\ang{s,t}\in P}.$ 
Obviously $(\ism\cG)\eqb(\ihf),$ thus, we have to prove that 
${\ism\cL}\reb{\ihf}$ for any $\cL$.

An action $\circ$ of $\isg$ on $\hfn$ 
is defined as follows: 
$g\circ n=g(n)$ for any $n\in\dN,$ and, by \dd\ve induction,  
$f\circ\ans{a_1,...,a_n})=\ans{f\circ a_1,...,f\circ a_n}$ 
for all $a_1,...,a_n\in\hfn.$ 
If $g\in\isg$ then $a\longmapsto g\circ a$ is a  
\dd\ve isomorfism of $\hfn$.

\ble
\label{grafs1}
Suppose that\/ $X,\,Y\sq\hfn$ are\/ \dd\ve transitive 
subsets of\/ $\hfn,$ the sets\/ $\dN\dif X$ and\/ 
$\dN\dif Y$ are infinite, and\/ 
${\ve\res X}\ihf{\ve\res Y}.$ 
Then there is a permutation\/ $f\in\isg$ such that\/ 
$Y=f\circ X=\ans{f\circ s:s\in X}$.
\ele 
\bpf
It follows from the assumption 
${\ve\res X}\ism\hfn{\ve\res Y}$ 
that there is an \dd\ve iso\-mor\-phism $\pi:X\onto Y.$  
Easily $\pi\res(X\cap\dN)$ is a bijection of $X_0=X\cap\dN$ 
onto $Y_0=Y\cap\dN,$ hence, there is $f\in\isg$ such that 
$f\res X_0=\pi\res X_0,$ and then we have $f\circ s=\pi(s)$ 
for any $s\in X$.
\epF{Lemma}

Coming back to the proof of Theorem~\ref{grafs}, we first 
show that ${\ism{\cG(m)}}\reb{\ihf}$ for any $m\ge 3,$ where 
$\cG(m)$ is the language with a single \dd mary predicate.  
We observe that $\ang{i_1,...,i_m}\in\hfn$ whenever  
$i_1,...,i_m\in\dN.$ 
Put $\vT(x)=\ans{\vt(s):s\in x}$ for every 
$x\in\mox{\cG(m)}=\cP(\dN^m),$ where 
$\vt(s)=\tce{\ans{\ang{2i_1,...,2i_m}}}$ for each 
$s=\ang{i_1,...,i_m}\in \dN^m,$ and finally, for $X\sq\hfn,$ 
$\tce X$ is the least \dd\ve transitive set $T\sq\hfn$ with 
$X\sq T.$  
It easily follows from Lemma~\ref{grafs} that 
$x\ism{\cG(m)}y$ is equivalent to  
${\ve\res\vT(x)}\ihf{\ve\res\vT(y)},$ which ends the 
proof of ${\ism{\cG(m)}}\reb{\ihf}$. 

It remains to show that ${\ism{\cL'}}\reb{\ihf},$ where 
$\cL'$ is the language with infinitely many binary 
predicates. 
In this case $\mox{\cL'}=\cP(\dN^2)^\dN,$ so that we can 
assume that every $x\in\mox{\cL'}$ has the form 
$x=\sis{x_n}{n\ge1},$ with 
$x_n\sq(\dN\dif\ans0)^2$ for all $n.$ 
Let 
$\vT(x)=\ans{s_n(k,l):n\ge1\land\ang{k,l}\in x_n}$ for 
any such $x,$ where 
\dm
s_n(k,l)=\tce{
\ans{\{...\{\ang{k,l}\}...\}\,,\,0}}
\,,\,\text{ with }\,n+2\,\text{ pairs of brackets }\,
\{\,,\,\}\,.
\dm
Then $\vT$ is a continuous reduction of $\ism{\cL'}$ 
to $\ihf$.
\epF{Theorem}

\punk{Proof of Theorem~\ref H}
\label{2t}

The proof (a version of the proof in \cite{frid}) is based on 
Scott's analysis. 

Define a family of Borel binary relations $\rrt\al st$ on 
$\cP(\dN^2),$ where $\al<\omi$ and $s,\,t\in\dN\lom,$ as 
follows:
\bit
\item\msur
$\rrq0stAB$ \ iff \ 
$A(s_i,s_j)\eqv B(t_i,t_j)$ for all $i,\,j<\lh s=\lh t$;

\item\msur
$\rrq{\al+1}stAB$ \ iff \ 
$\kaz k\:\sus l\:(\rrQ\al{s\we k}{t\we l}AB)$ 
and $\kaz l\:\sus k\:(\rrQ\al{s\we k}{t\we l}AB)$;

\item 
if $\la<\omi$ is limit then: \ 
$\rrq\la stAB$ \ iff \ $\rrq\al stAB$ for all $\al<\la$.
\eit
We define $\rro\al stAB$ iff $\rrq\al stAB\,;$ then, by 
induction on $\al,$ each $\rrO\al$ is easily a Borel \er\ 
on $\dN\lom\ti\cP(\dN^2),$ and   
${\rrO\ba}\sq{\rrO\al}$ whenever $\al<\ba$. 

Let $\rE=\ergx$ be a \er\ induced by a Polish action of 
a Polish group $\dG$ on a Polish space $\dX.$ 
Suppose that $\rE$ is reducible to a Polish action of 
$\isg$ by a \bm\ map. 
According to Theorems \ref{siy2cc}, \ref{grafs}, and 
Proposition~\ref{BBM}, 
there is a \bm\ reduction $\vt:\dX\to\cP(\dN^2)$ of $\rE$ to 
$\ism\cG.$ 
The reduction is continuous on a dense $\Gd$ set $D_0\sq\dX.$ 
Recall that, for $A,\,B\sq\dN^2,$ $A\ism\cG B$ means that 
there is $f\in\isg$ with $A(k,l)\eqv B(f(k),f(l))$ for all 
$k,\,l.$  
We easily prove ${\ism\cG}\sq{\rrt\al st},$ 
where $t=f\circ s,$ by induction on $\al,$ 
in particular, ${\ism\cG}\sq{\rrt\al\La\La},$
where $\La$ is the empty sequence. 
Since $\vt$ is a reduction, the equivalence  
${x\rE y}\eqv{\vt(x)\ism\cG\vt(y)}$ holds for all $x,\,y.$ 
Our goal is to find a $\Gd$ dense set $D\sq D_0$ 
and an ordinal $\al<\omi$ such that
\bit
\item[$\mtho(\ast)$] 
the implication 
${x\nE y}\imp{\nrq\al\La\La{\vt(x)}{\vt(y)}}$ \ 
holds for all $x,\,y\in D$.
\eit

To find $D$ fix a countable transitive model $\mm$ of 
$\zhc$ (see above). 
We assume that $\dX,$ the group $\dG,$ the action, 
$D_0,\msur$ $\vt\res D_0$ are assumed to be coded in $\mm$ in 
the same sense as in the proof of Lemma~\ref{L1}. 
We assert that the set $D$ of all Cohen generic, over $\mm,$ 
points of $\dX$ 
(a dense $\Gd$ subset of $\dX$ included in $D_0$) 
satisfies $(\ast)$.

Suppose that $x,\,y\in D.$ 
First consider the case when $\ang{x,y}$ 
is a Cohen generic pair over $\mm.$  
If ${x\nE y}$  then, by the choice of $\vt,$ we have 
$\vt(x)\not\ism\cG \vt(y),$ hence, 
this fact holds in $\mm[x,y]$ by the Mostowski 
absoluteness. 
Therefore, arguing in $\mm[x,y]$ 
(which is still a model of $\zhc,$ see Footnote~\ref{mxy}), 
we find an ordinal 
$\al\in\Ord^\mm=\Ord^{\mm[x,y]}$ with 
$\nrq\al\La\La {\vt(x)}{\vt(y)}.$ 
Moreover, since the Cohen forcing satisfies \ccc, 
there is an ordinal $\al\in\mm$ such that 
$\nrq\al\La\La {\vt(x)}{\vt(y)}$ for \poq{all} 
Cohen generic, over $\mm,$ pairs $\ang{x,y}\in D^2$ 
with $x\nE y.$   
It remains to show that this also holds when $x,\,y\in D$ 
with $x\nE y$ do not form a Cohen generic pair. 

Let $g\in\dG$ be Cohen generic over $\mm[x,y].$
Then $z=g\app x\in\dX$ is easily Cohen generic over 
$\mm[x,y]$ (because the action is continuous), furthermore, 
$x\rE z,$ hence, $y\nE z.$  
However $y$ is generic over $\mm$ and $z$ is generic over 
$\mm[y],$ thus, $\ang{y,z}$ is Cohen generic over $\mm,$ 
hence, we 
have $\nrq\al\La\La {\vt(z)}{\vt(y)}$ by the above. 
On the other hand, $\rrq\al\La\La {\vt(x)}{\vt(z)}$ holds 
because $x\rE z,$ thus, we finally obtain 
$\nrq\al\La\La {\vt(x)}{\vt(y)},$ as required by $(\ast)$.

To conclude, we have  
${x\rE y}\eqv{\rrq\al\La\La{\vt(x)}{\vt(y)}}$   
for all $x,\,y\in D.$ 
In this case we can easily redefine $\vt$ on the  
complement of $D$ in $\dX$ so that the equivalence holds 
for all $x,\,y\in\dX,$ in other words, the improved $\vt$ 
is a BM (because $\vt\res D$ is 
continuous and $D$ is a dense $\Gd$) 
reduction of $\rE$ to $\rrt\al \La\La.$ 

The following result completes the proof of the theorem.

\bpro
\label{r2t}
Any \er\/ $\rrO\al$ is Borel reducible to some\/ $\rT_\ga$.
\epro
\bpf 
We have ${\rrO0}\reb{\rT_0}$ since $\rrO0$ has countably 
many equivalence classes, all of which are 
open-and-closed sets. 
To carry out the step $\al\mapsto \al+1$ note that the map   
$\ang{s,A}\mapsto \sis{\ang{s\we k,A}}{k\in\dN}$ is a Borel 
reduction of $\rrO{\al+1}$ to $(\rrO\al)^\iy.$ 
As for the limit step, let 
$\la=\ans{\al_n:n\in\dN}$ be a limit ordinal, and 
${\rR}=\bigvee_{n\in\dN}{\rrO{\al_n}},$ \ie, $\rR$ is 
a \er\ on $\dN\ti\dN\lom\ti\cP(\dN^2)$ defined so that 
$\ang{m,s,A}\rR\ang{n,t,B}$ iff $m=n$ and $\rrq{\al_m} stAB.$ 
However the map $\ang{s,A}\mapsto\sis{\ang{m,s,A}}{m\in\dN}$ 
is a Borel reduction of $\rrO\la$ to $\rR^\iy.$
\epf

\qeDD{Theorem~\ref{H}}

\parf{\PPP\ \er s and irreducibility of $\rtd$}
\label{rtd}

This section contains a theorem saying that the \er\ $\rtd$ 
of equality of countable sets of the reals is not Borel 
reducible to \er s which belong to a family of 
{\it \PP\/} \er s, including, for 
instance, continuous actions of \cli\ groups, some 
ideals, not only Polishable, together with \er s having 
$\Gds$ equivalence classes, 
and is closed under the Fubini product modulo $\ifi.$   
The definition of the family is based 
on a rather metamathematical property which we extracted 
from Hjorth~\cite{h:orb}.

\punk{\PPP\ \er s}
\label{ppdef}

First of all, if $X$ is an analytic set in the universe 
$\dV$ of all sets 
(in particular, this applies when $X$ is Borel), 
and $\dvp$ is a generic extension $\dV,$   
then $\di X$ will denote the result of the sequence of 
operations contained in the definition of $X$ but applied 
in $\dvp.$ 
The correctness of this definition follows 
from the Shoenfield absoluteness 
theorem, and easily $X=\di X\cap\dV$. 

For instance, if, in $\dV,$ $\rE$ is an analytic \er\ on a 
polish space $\dX,$ then, still by the Shoenfield 
absoluteness, $\drE$ is an analytic \er\ on $\di\dX.$  
If now $x\in \dX$ (hence, $x\in\dV$) then the \dde class 
$\ek x\rE\sq\dX$ of $x$(defined in $\dV$) is included in 
a unique \dd\drE class $\ek x{\drE}\sq\di\dX$ (in $\dvp$).  
Classes of the form $\ek x{\drE}\zT x\in\dX,$ belong to a 
wider category of \dd\drE classes which admit a description 
from the \dd\dV th point of view. 

\bdf
\label{virt}
Assume that $\dP$ is a notion of forcing in $\dV.$ 
A {\it virtual\/ \dde class\/} is any \ddp term $\xib$ 
such that $\dP$ forces $\xib\in\di\dX$ and 
$\dPP$ forces $\tal\drE \tar$.~\footnote
{\label{xilr}\ $\tal$ and $\tar$ are \dd\dPP terms  
meaning $\xib$ associated with the resp.\ left and right 
factors $\dP$ in the product forcing, formally, 
$\tal[U\ti V]=\xib[U]$ and $\tar[U\ti V]=\xib[V]$ for any 
\dd\dPP generic set $U\ti V,$ where, say, $\xib[U]$ is the 
interpretation of a term $\xib$ via a generic set $U$.}  
A virtual class is {\it\PP\/} if there is, in $\dV,$ a 
point $x\in\dX$ which pins it in the sense that 
$\dP$ forces $x\drE \xib.$  
Finally, an analytic \er\ $\rE$ is {\it\PP\/} if, 
for any forcing notion $\dP\in\dV,$ 
all virtual \dde classes are \PP.
\edf

If $\xib$ is a virtual \dde class then, in any extension 
$\dvp$ of $\dV,$ if $U$ and 
$V$ are generic subsets of $\dP$ then $x=\xib[U]$ 
and $y=\xib[V]$ belong to $\di\dX$ and satisfy $x\drE y,$ 
hence, $\xib$ induces a \dd\drE class in the 
extension. 
If $\xib$ is \PP\ then this class contains an element in the 
ground universe $\dV$ --- in other words, \PP\ virtual 
classes induce \dd\drE equivalence classes of the form 
$\ek x\drE\zT x\in\dV$ in the extensions of the universe 
$\dV.$ 

We prove below that $\rtd$ is not \PP, moreover, $\rtd$ is 
not Borel reducible to any \PP\ analytic \er. 
In addition, we give a simplified proof of Hjorth's 
theorem that continuous actions of Polish \cli\ groups never 
induce \PP\ orbit \er s, introduce a family of \PP\ \er s 
associated with $\Fsd$ ideals, show that any Borel \er\ 
whose all equivalence classes are $\Gds$ is \PP, 
and prove that 
the class of all \PP\ analytic \er s is closed under the 
Fubini product over $\ifi$.

\punk{\PPP\ \er s do not reduce $\rtd$}
\label{tagcli}

Recall that, modulo $\eqb\,,\msur$ $\rtd$ is a \er\   
on $\dNN$ defined as follows: 
$x\rtd y$ iff $\ran x=\ran y.$

\ble
\label{t2unlike} 
$\rtd$ is \poq{not} \PP.
If\/ $\rE,\,\rF$ are analytic\/ \er s, $\rE\reb\rF,$  
and\/ $\rF$ is \PP, then so is\/ $\rE.$ 
Hence,\/ $\rtd$ is not Borel reducible to a \PP\ 
analytic \er. 
\ele
\bpf
To prove that $\rtd$ is not \PP, consider 
$\dP=\coll(\dN,\dn),$ a forcing to produce a generic 
map $f:\dN\onto\dn.$ 
($\dP$ consists of all functions $p:u\to\dn$ where  
$u\sq\dN$ is finite.)
Let $\xib$ be a \ddp term for the set 
$\ran f=\ans{f(n):n\in\dN}.$ 
Then $\xib$ is obviously a virtual \dd\rtd class, but 
it is not \PP\ because $\dnn$ is 
uncountable in the ground universe $\dV$.  

Suppose that, in $\dV,$ $\vt:\dX\to\dY$ is a Borel
reduction of $\rE$ to $\rF,$ where $\dX=\dom\rE$ and 
$\dY=\dom\rF.$    
We can assume that $\dX$ and $\dY$ are just two copies 
of $\dn.$ 
Let $\dP$ be a forcing notion and a \ddp term $\xib$ be a 
virtual \dde class.  
By the Shoenfield absoluteness, $\di\vt$ is a 
reduction of $\di\rE$ to $\di\rF$ in any extension of $\dV,$ 
hence, $\sgb,$ a \ddp name for $\di\vt(\xib),$ is also a 
virtual \ddf class. 
Since $\rF$ is \PP, there is $y\in\dY$ such that 
$\dP$ forces $y\drF\sgb.$ 
Note that it is true in the \ddp extension that 
$y\drF \di\vt(x)$ for some $x\in\di\dX,$ hence, by Shoenfield, 
in the ground universe there is $x\in\dX$ with $y\rF\vt(x).$ 
Clearly $\dP$ forces $x \drE \xib$.
\epf

\punk{Fubini product of \PP\ \er s is \PP}
\label{fbbsbs}

Recall that the Fubini product $\rE=\fps{\ifi}{\rE_k}{k\in\dN}$ 
of \er s $\rE_k$ on $\dX_k$ modulo $\ifi$ is a \er\ on 
$\dX=\prod_k\dX_k$ defined as follows: 
$x\rE y$ if $x(k)\rE_k y(k)$ for all but finite $k$.

\ble
\label{bsfub}
The family of all analytic \PP\/ \er s is closed 
under Fubini products modulo\/ $\ifi$.
\ele
\bpf
Suppose that analytic \er s $\rE_k$ on Polish spaces 
$\dX_k$ are \PP; prove that the Fubini product 
$\rE=\fps{\ifi}{\rE_k}{k\in\dN}$ is a \PP\ \er\ on 
$\dX=\prod_k\dX_k.$ 
Consider a forcing notion $\dP$ and a \ddp term $\xib$ 
which is a virtual \dde class.  
There is a number $k_0$ and conditions $p,\,q\in \dP$ 
such that $\ang{p,q}$ \dd\dPP forces 
$\tal(k)\mathbin{\di{\rE_k}}\tar(k)$ for all $k\ge k_0.$ 
As all $\rE_k$ are \er s, we conclude that the condition  
$\ang{p,p}$ also forces 
$\tal(k)\mathbin{\di{\rE_k}}\tar(k)$ for all $k\ge k_0.$
Therefore, since $\rE_k$ are \PP, there exists, in $\dV,$ 
a sequence of points $x_k\in\dX_k$ such that $p$ \ddp forces 
$x_k \mathbin{\di{\rE_k}}\xib(k)$ for any $k\ge k_0.$ 
Let $x\in\dX$ satisfy $x(k)=x_k$ for all $k\ge k_0.$ 
(The values $x(k)\in\dX_k$ for $k<k_0$ can be arbitrary.)  
Then $p$ obviously \ddp forces $x\drE \xib.$ 

It remains to show that just every $q\in\dP$ also forces 
$x\drE \xib.$ 
Suppose otherwise, \ie, some $q\in\dP$ forces that 
$x\drE \xib$ \poq{fails}. 
Consider the pair $\ang{p,q}$ as a condition in $\dPP:$ 
it forces $x\drE\tal$ and $\neg\;{x\drE\tar},$ as 
well as $\tal\drE\tar$ by the choice of $\rE$ and 
$\xib,$ which is a contradiction.
\epf

\punk{Left-invariant actions induce \PP\ \er s}
\label{bs:cli}

Recall that a Polish group $\dG$ is 
{\it complete left-invariant\/}, \cli\ for brevity, if 
$\dG$ admits a compatible left-invariant complete metric. 
Then easily $\dG$ also admits a compatible 
\poq{ri}g\poq{ht}-invariant complete metric, 
which will be practically used.

\bte[{{\rm Hjorth~\cite{h:orb}}}]
\label{clit}
Any \er\/ $\rE=\egx$ induced by a Polish action of 
a\/ \cli\ group\/ $\dG$ on a Polish space\/ $\dX$  
is \PP, hence, 
$\rtd$ is not Borel reducible to\/ $\rE$. 
\ete
\bpf
Let $\dP$ be a forcing notion and $\xib$ be a virtual 
\dde class. 
Let $\le$ denote the partial order of $\dP\,;$ we assume 
that $p\le q$ means that $p$ is a stronger condition.
Let us fix a compatible complete right-invariant 
metric $\rho$ on $\dG.$ 
For any $\ve>0,$ put $G_\ve=\ans{g\in\dG:\rho(g,1_\dG)<\ve}.$ 
Say  that $q\in\dP$ {\it is of size\/ $\le\ve$\/} 
if $\ang{q,q}$ \dd{\dPP}forces the existence of 
$g\in \di{G_\ve}$ such that $\tal=g\app \tar$. 

\ble
\label{l4}
If\/ $q\in\dP$ and\/ $\ve>0,$ then there exists a condition\/ 
$r\in\dP\zT r\le q,$ of size\/ $\le \ve$.
\ele
\bpf
Otherwise for any $r\in\dP\zT r\le q,$ there 
is a pair of conditions $r',\,r''\in\dP$ stronger than 
$r$ and such that $\ang{r',r''}$ \dd{\dPP}forces that 
there is no $g\in \di{G_\ve}$ with $\tal=g\app\tar.$ 
Applying an ordinary splitting construction in such 
a generic extension $\dvp$ of $\dV$ where 
$\cP(\dP)\cap\dV$ is countable, we find an uncountable 
set $\cU$ of generic sets $U\sq\dP$ with $q\in U$ such 
that any pair $\ang{U,V}$ of $U\ne V$ in $\cU$ is 
\dd\dPP generic (over $\dV$), hence, there is no 
$g\in \di{G_\ve}$ with 
$\xib[U]=g\app \xib[V].$~\footnote
{\label{tauu}\ $\xib[U]$ is the interpretation of the 
\ddp term $\xib$ obtained by taking $U$ as the generic 
set.} 
Fix $U_0\in \cU.$
We can associate, in $\dvp,$ with each $U\in \cU,$ an 
element $g_U\in\di G$ such that 
$\xib[U]=g_U\app \xib[U_0];$ then $g_U\nin\di{G_\ve}$ 
by the above. 
Moreover, we have $g_Vg_U\obr\app \xib[U]=\xib[V]$ 
for all $U,\,V\in \cU,$ 
hence $g_Vg_U\obr\nin\di{G_\ve}$ whenever $U\ne V,$ which 
implies $\rho(g_U,g_V)\ge\ve$ by the right invariance. 
But this contradicts the separability of $G$.
\epF{Lemma}

Coming back to the theorem, suppose on the contrary that 
a condition $p\in\dP$ forces that there is no $x\in\dX$ 
(in the ground universe $\dV$) 
satisfying $x\drE\xib.$ 
According to the lemma, there is, in $\dV,$ a sequence 
of conditions $p_n\in\dP$ of size $\le 2^{-n},$ 
and closed sets $X_n\sq\dX$ with \dd\dX diameter 
$\le 2^{-n},$ such that 
$p_0\le p\zT\,p_{n+1}\le p_n\zT\,X_{n+1}\sq X_n,$ 
and $p_n$ forces $\xib\in\di{X_n}$ for any $n.$ 
Let $x$ be the common point of the sets $X_n$ in $\dV.$ 
We claim that $p_0$ forces $x\drE \xib$. 

Indeed, otherwise there is $q\in\dP\zT q\le p_0,$ which 
forces $\neg\;{x\drE \xib}.$
Consider an extension $\dvp$ of $\dV$ rich enough to 
contain, for any $n,$ a generic set $U_n\sq\dP$ with 
$p_n\in U_n$ such that each pair $\ang{U_n,U_{n+1}}$ 
is \dd\dPP generic 
(over $\dV$), and, in addition, $q\in U_0.$    
Let $x_n=\xib[U_n]$ (an element of $\di\dX$). 
Then $\sis{x_n}\to x.$ 
Moreover, for any $n,$ both $U_n$ and $U_{n+1}$ 
contain $p_n,$ hence, as $p_n$ has size 
$\le 2^{-n-1},$ there is $g_{n+1}\in \di{\dG_\ve}$ with 
$x_{n+1}=g_{n+1}\app x_n.$ 
Thus, $x_n=h_n\app x_0,$ where $h_n=g_{n}...g_1.$ 
Note that 
$\rho(h_n,h_{n-1})=\rho(g_n,1_\dG)\le 2^{-n+1}$ by 
the right-invariance of the metric, thus, 
$\sis{h_n}{n\in\dN}$ is a Cauchy sequence in $\di\dG.$ 
Let $h=\tlim_{n\to\iy}h_n\in\di\dG$ be its limit. 
As the action is continuous, we have 
$x=\tlim_nx_n=h\app x_0.$ 
It follows that $x\drE x_0$ holds in $\dvp,$ hence, 
also in $\dV[U_0].$ 
However $x_0=\xib[U_0]$ while $q\in U_0$ forces 
$\neg\;{x\drE \xib},$ which is a contradiction. 

Thus $p_0$ \ddp forces $x\drE \xib.$ 
Then any $r\in\dP$ also forces $x\drE \xib:$ 
indeed, if some $r\in\dP$ forces $\neg\;{x\drE \xib}$ 
then the pair $\ang{p_0,r}$ forces, in $\dPP,$ that 
$x\drE \tal$ and $\neg\;{x\drE \tar},$ which 
contradicts the fact that $\dPP$ forces 
$\tal\drE\tar$.\vtm

\epF{Theorem~\ref{clit}}

\punk{All \er s with $\Gds$ classes are \PP}
\label{s02}

We have a non-\PP\ \er\ $\rtd,$ obviously of class $\Fsd;$ 
the following theorem shows that this is the simplest 
possible case of non-\PP\ \er s.

\bte
\label{fspp}
Any Borel \er\/ $\rE$ whose all equivalence classes are\/ 
$\Gds$ is \PP.
\ete
\bpf[{{\rm Based on an idea communicated by Hjorth}}]
We can assume that $\dom\rE=\dnn.$ 
It follows from a theorem of Louveau~\cite{louv}, that 
there is a Borel map $\ga,$ defined on $\dnn$ so that 
$\ga(x)$ is a \dd\Gds code of $\ek x\rE$ for any 
$x\in\dnn,$ that is, for instance, 
$\ga(x)\sq\dN^2\ti\dN\lom$ and 
$\ek x\rE=
\bigcup_i\bigcap_j\bigcup_{\ang{i,j,s}\in\ga(x)}B_s,$ 
where $B_s=\ans{a\in\dnn:s\su a}$ for all $s\in\dN\lom.$  

Let $\dP=\stk\dP\le$ be a forcing notion, 
and $\xib$ be a virtual \dde class, thus, $\dPP$ 
forces $\tal\drE\tar,$ hence, there is a number $i_0$ 
and a condition $\ang{p_0,q_0}\in\dPP$ which forces 
$\tal\in\di\vt(\tar),$ where 
$\vt(x)=\bigcap_j\bigcup_{\ang{i_0,j,s}\in\ga(x)}B_s$ for 
all $x\in\dnn.$ 

The key idea of the proof is to substitute $\dP$ by the 
Cohen forcing. 
Let $\dS$ denote the set of all $s\in\dN\lom$ such that 
$p_0$ does \poq{not} \ddp force that $s\not\su \xib.$ 
We consider $\dS$ as a forcing, and $s\sq t$ 
(\ie, $t$ is an extension of $s$) means that $t$ is a 
stronger condition; $\La,$ the empty sequence, 
is the weakest condition in $\dS.$  
If $s\in \dS$ then obviously there is at least one $n$ 
such that $s\we n\in \dS,$ hence, $\dS$ forces an element 
of $\dnn,$ whose \dd\dS name will be $\fa$. 

\ble
\label{s0}
The pair\/  $\ang{\La,q_0}$ \dd{\dSP}forces\/ 
$\fa\in\di\vt(\xib).$ 
\ele
\bpf
Otherwise some condition $\ang{s_0,q}\in\dSP$ with  
$q\le q_0$ forces $\fa\nin\di\vt(\xib).$ 
By the definition of $\vt$ we can assume that there 
is $j_0$ satisfying  
\dm
\ang{s_0,q}
\quad \text{\dd{\dSP}forces}\quad
\neg\;\sus s\:\skl
\ang{i_0,j_0,s}\in\ga(\xib)\land s\su\fa
\skp.
\eqno(\ast)
\dm
Since $s_0\in\dS,$ there is a condition 
$p'\in\dP\zq p'\le p_0,$ which \ddp forces $s_0\su\xib.$ 
By the choice of $\ang{p_0,q_0}$ we can assume that, 
for some $s\in\dS$ and $q'\in\dP\zq q'\le q$,
\dm
\ang{p',q'}
\quad \text{\dd{\dPP}forces}\quad
\ang{i_0,j_0,s}\in\ga(\tar)\land s\su\tal\,.
\dm
This means that 
$1)\msur$ $p'$ \ddp forces $s\su\xib$ and 
$2)\msur$ $q'$ \ddp forces $\ang{i_0,j_0,s}\in\ga(\xib).$ 
In particular, by the above, $p'$ forces both $s_0\su\xib$ 
and $s\su\xib,$ therefore, 
either $s\sq s_0$ -- then let $s'=s_0,$ 
or $s_0\su s$ -- then let $s'=s.$ 
In both cases, $\ang{s',q'}$ \dd\dSP forces 
$\ang{i_0,j_0,s}\in\ga(\xib)$ and $s\su\fa,$ 
contradiction to $(\ast)$.
\epF{Lemma} 

Note that $\dS$ is a subforcing of the Cohen forcing 
$\dC=\dN\lom,$ therefore, by the lemma, there is a 
\dd\dC term $\sgb$ such that $\ang{\La,q_0}$ \dd\dCP 
forces 
$\sgb\in\di\vt(\xib),$ hence, forces $\sgb\drE\xib.$ 
It follows, by consideration of the forcing $\dCP\ti\dP,$ 
that generally $\dCP$ forces $\sgb\drE\xib.$ 
Therefore, by ordinary arguments, first,  
$\dCC$ forces $\sgbl\drE\sgbr,$ and second, to prove the 
theorem it suffices now to find $x\in\dnn$ in $\dV$ such 
that $\dC$ forces $x\drE\sgb.$ 
This is our next goal.

Let $\fa$ be the \dd\dC name of the Cohen generic element 
of $\dnn.$ 
The term $\sgb$ can be of arbitrary nature, but we can 
substitute it by a term of the form $\di f(\fa),$ where 
$f:\dnn\to\dnn$ is a Borel map in the ground universe $\dV.$ 
It follows from the above that 
$\di f(\fa)\drE \di f(\fb)$ for any \dd\dCC generic, 
over $\dV,$ pair $\ang{\fa,\fb}\in\dnn\ti\dnn.$ 
We conclude that $\di f(\fa)\drE \di f(\fb)$ 
also holds even for any pair 
of separately Cohen generic $\fa,\,\fb\in\dnn.$ 
Thus, in a generic extension of $\dV,$ where there exist 
comeager-many Cohen generic reals, there is a comeager 
$\Gd$ set $X\sq\dnn$ such that $\di f(a)\drE \di f(b)$ 
for all $a,\,b\in X.$ 
By the Shoelfield absoluteness, the statement of existence 
of such a set $X$ is true also in $\dV,$ hence, in $\dV,$ 
there is $x\in\dnn$ such that we have $x\rE f(a)$  
for comeager-many $a\in\dnn.$ 
This is again a Shoenfield absolute property of $x,$ 
hence, $\dC$ forces $x\drE \di f(\fa),$ as required.\vtm

\epF{Theorem~\ref{fspp}}

\punk{A family of \PP\ ideals}
\label{bsti}

Let us say that a Borel ideal $\cI$ is {\it\PP\/} 
if so is the induced \er\ $\rei.$ 
It follows from Theorem~\ref{clit} that any 
P-ideal is \PP\ because Borel P-ideals are polishable  
\cite{sol'} while all Polish Abelian groups are \cli. 
Yet there exist non-P \PP\ ideals. 

We introduce here a family of such ideals. 
Suppose that $\sis{\vpi_i}{i\in\dN}$ is a sequence of 
lower semicontinuous (\lsc) submeasures on $\dN.$ 
Define
\dm
\Exh_{\sis{\vpi_i}{}}\,=\,
\ans{X\sq\dN:\vpy(X)=0}\,,
\quad\hbox{where}\quad
\vpy(X)=\tlis_{i\to\iy}\vpi_i(X)\,.
\dm
the exhaustive ideal of the sequence of submeasures. 
By Solecki's Theorem~\cite{sol'} for any Borel P-ideal 
$\cI$ there is an \lsc\ submeasure $\vpi$ such that  
$\cI=\Exh_{\sis{\vpi_i}{}}=\Exh_\vpi,$ where 
$\vpi_i(x)=\vpi(x\cap\iry i),$ however, for example, the 
non-polishable ideal $\Ii=\fio$ also is of the form 
$\Exh_{\sis{\vpi_i}{}},$ where for $x\sq\dN^2$ we define 
$\vpi_i(x)=0\,\hbox{ or }\,1$ if resp.\ 
$x\sq\,\hbox{ or }\,\not\sq\ans{0,...,n-1}\ti\dN$. 

\bte
\label{ibs}
All ideals of the form\/ $\Exh_{\sis{\vpi_i}{}}$ are 
\PP.
\ete
\bpf
Let $\cI=\Exh_{\sis{\vpi_i}{}},$ where all $\vpi_i$ 
are \lsc\ submeasures on $\dN.$  
We can assume that the submeasures $\vpi_i$ decrease, 
\ie, $\vpi_{i+1}(x)\le\vpi_i(x)$ for any $x,$ for if not 
consider the \lsc\ submeasures 
$\vpi'_i(x)=\tsup_{j\ge i}\vpi_j(x).$ 

Suppose that $\rE=\rei$ is not \PP.
Then there is a forsing notion $\dP,$ a virtual 
\dde klass $\xib,$ and a condition $p\in \dP$ which 
\ddp forces $\neg\;x\drE\xib$ for any $x\in\pn$ in 
$\dV.$ 
By definition, for any $p'\in\dP$ i $n\in\dN$ there 
are $i\ge n$ and conditions $q,\,r\in\dP$ with 
$q,\,r\le p',$ such that $\ang{q,r}$ \dd{\dPP}forces  
the inequality $\vpi_i(\tal\sd\tar)\le 2^{-n-1},$ 
hence, $\ang{q,q}$ \dd{\dPP}forces  
$\vpi_i(\tal\sd\tar)\le 2^{-n}.$ 
It follows that, in $\dV,$ there is a sequence of numbers 
$i_0<i_1<i_2<...,$ and a sequence 
$p_0\ge p_1\ge p_2\ge...$ of conditions in $\dP,$ and, 
for any $n,$ a set $u_n\sq\ir0n,$ such that $p_0\le p$ i 
\ben
\tenu{(\arabic{enumi})}
\itla{bs1}
each $p_n$ \dd\dP forces $\xib\cap\ir0n=u_n$;

\itla{bs2}
each $\ang{p_n,p_n}$ \dd\dPP forces  
$\vpi_{i_n}(\tal\sd \tar)\le 2^{-n}.$ 
\een
Arguing in $\dV,$ put $a=\bigcup_nu_n;$ then  
$a\cap\ir0n=u_n$ for all $n.$ 
We claim that $p_0$ forces $a\drE\xib,$ contrary to the 
assumption above, which proves the theorem. 

Indeed, otherwise there is a condition $q_0\le p_0$ 
which forces $\neg\;a\drE\xib.$ 
Consider a generic extension $\dvp$ of the universe, where 
there is a sequence of \ddp generic sets $U_n\sq\dP$ such 
that, for any $n,$ the pair $\ang{U_n,U_{n+1}}$ is 
\dd\dPP generic, $p_n\in U_n,$ and in addition 
$q_0\in U_0.$ 
Then, in $\dvp,$ the sets $x_n=\xib[U_n]\in\pn$ satisfy
$\vpi_{i_n}(x_n\sd x_m)\le 2^{-n}$ by \ref{bs2}, 
whenever $n\le m.$ 
It follows that $\vpi_{i_n}(x_n\sd a)\le 2^{-n},$ because 
$a=\tlim_mx_m$ by \ref{bs1}. 
However we assume that the submeasures $\vpi_j$ decrease, 
hence, $\vpy(x_n\sd a)\le 2^{-n}.$ 
On the other hand, $\vpy(x_n\sd x_0)=0$ because $\xib$ 
is a virtual \dde class.
We conclude that $\vpy(x_0\sd a)\le 2^{-n}$ 
for any $n,$ in other words, $\vpy(x_0\sd a)=0,$ that is,  
$x_0\drE a,$ which is a contradiction with the choice of 
$U_0$ because $x_0=\xib[U_0]$ and $q_0\in U_0$.
\epf

\subsubsection*{\bfsl Questions}

\bqe
\label{ppq1}
Are all Borel ideals \PP~? \
The expected answer ``yes'' would show that $\rtd$ is not 
Borel reducible to any Borel ideal.  
Moreover, is any orbit \er\ of a Borel action of a Borel 
\cli\ group \PP~? 
\eqe

\bqe[{{\rm Kechris}}]
\label{ppq2}
Is there a \dd\reb least non-\PP\ Borel \er~? \ 
It was once expected that $\rtd$ is such, but Hjorth 
informed us that there is a strictly \dd\reb smaller 
non-\PP\ Borel \er\ of a rather complicated nature.  
\eqe

\subsubsection*{\bfsl Acknowledgements}

We are thankful to Greg Hjorth, A.~S.~Kechris, and 
Su Gao for useful discussions related to the content 
of this paper. 
We also are thankful to Greg Hjorth for his kind 
permission to include Theorem~\ref{fspp} in this 
paper.

\end{document}